\documentclass[11pt,a4paper]{article}
\usepackage{amsmath,amssymb}
\usepackage{latexsym}
\usepackage{mathrsfs}
\usepackage{amsfonts}
\usepackage{hyperref}
\usepackage{amsthm}
\usepackage{appendix}
\usepackage[square,comma,sort&compress,numbers]{natbib}
\usepackage{color}

\hfuzz=\maxdimen
\theoremstyle{definition}
\newtheorem{lemma}{Lemma}[section]
\newtheorem{definition}[lemma]{Definition}
\newtheorem{theorem}[lemma]{Theorem}
\newtheorem{corollary}[lemma]{Corollary}
\newtheorem{proposition}[lemma]{Proposition}
\newtheorem{remark}{Remark}

\DeclareFixedFont{\Acknowledgment}{OT1}{cmr}{bx}{n}{14pt}
\begin{document}

\title{\bf Parameterized discrete uniformization theorems and curvature flows for polyhedral surfaces, II}
\author{Xu Xu, Chao Zheng}
\maketitle

\begin{abstract}
This paper investigates the combinatorial $\alpha$-curvature for vertex scaling of
piecewise hyperbolic metrics on polyhedral surfaces,
which is a parameterized generalization of the classical combinatorial curvature.
A discrete uniformization theorem for combinatorial $\alpha$-curvature is established,
which generalizes Gu-Guo-Luo-Sun-Wu's discrete uniformization theorem for classical combinatorial curvature \cite{Gu2}.
We further introduce combinatorial $\alpha$-Yamabe flow and combinatorial $\alpha$-Calabi flow
for vertex scaling to find piecewise hyperbolic metrics with prescribed combinatorial $\alpha$-curvatures.
To handle the potential singularities along the combinatorial curvature flows, we do surgery along the flows by edge flipping.
Using the discrete conformal theory established by Gu-Guo-Luo-Sun-Wu \cite{Gu2},
we prove the longtime existence and convergence of combinatorial $\alpha$-Yamabe flow
and combinatorial $\alpha$-Calabi flow with surgery,
which provide effective algorithms for finding piecewise hyperbolic metrics with prescribed combinatorial $\alpha$-curvatures.
\end{abstract}

\textbf{Keywords}: Discrete uniformization theorem; Piecewise hyperbolic metric;
Combinatorial curvature; Combinatorial Yamabe flow; Combinatorial Calabi flow

\section{Introduction}
This is a continuation of \cite{Xu1}, in which the first author
established a parameterized discrete uniformization theorem for combinatorial $\alpha$-curvature for vertex scaling of piecewise linear metrics (PL metrics for short in the following)
and studied the longtime behavior of combinatorial $\alpha$-Yamabe flow and combinatorial $\alpha$-Calabi flow
for vertex scaling of PL metrics by surgery.
In this paper, we establish the counterpart results for combinatorial $\alpha$-curvature
for vertex scaling of piecewise hyperbolic metrics (PH metrics for short in the following) on surfaces.

%
%

Suppose $S$ is a connected closed surface and $V$ is a finite subset of $S$ with $|V|=N$, $(S, V)$ is called a marked surface.
A PH metric on the marked surface $(S,V)$ is a hyperbolic metric on $S$ with conic singularities contained in $V$.
Suppose $\mathcal{T}={(V,E,F)}$ is a triangulation of $(S, V)$, where $V,E,F$ represent the sets of vertices, edges and faces respectively,
$(S, V, \mathcal{T})$ is called a triangulated surface.
Denote a vertex, an edge and a face in the triangulation by $v_i,\{v_iv_j\},\triangle v_iv_jv_k$ respectively,
where ${i,j,k}$ are natural numbers.
A PH metric
on a triangulated surface $(S,V, \mathcal{T})$
defines a map $d: E\rightarrow (0, +\infty)$ such that for any triangle $\triangle v_iv_jv_k\in F$,
$d_{ij}, d_{ik}, d_{jk}$ satisfy the triangle inequalities, where $d_{st} : =d(\{v_sv_t\})$ for $\{v_sv_t\}\in E$.
Conversely, by gluing hyperbolic triangles isometrically along the edges in pair,
we obtain PH metrics on a triangulated surface $(S,V, \mathcal{T})$,
which give rise to PH metrics on the marked surface $(S, V)$.
The classical combinatorial curvature $K: V\rightarrow (-\infty, 2\pi)$ is used to describe the conic singularities of PH metrics at the vertices.
Suppose $(S, V, \mathcal{T})$ is a triangulated surface with a PH metric $d$,
the classical combinatorial curvature for $d$ at $v_i\in V$ is defined to be
\begin{equation}\label{K_i}
K_i=2\pi-\sum_{\triangle v_iv_jv_k\in F}\alpha^{jk}_i,
\end{equation}
where $\alpha^{jk}_i$ is the inner angle at the vertex $v_i$ in the hyperbolic triangle $\triangle v_iv_jv_k\in F$.
Note that the combinatorial curvature for the PH metric $d$ is independent of the geometric triangulations
of $d$ on $(S, V)$.

The vertex scaling for PL metric
on polyhedral surfaces is a discrete analogy of the conformal change in Riemannian geometry, which
was introduced physically by R\v{o}cek-Williams \cite{RW} and mathematically by Luo \cite{L1} independently.
Motivated by the vertex scaling for PL metrics on polyhedral surfaces,
Bobenko-Pinkall-Springborn \cite{BPS} introduced the following definition of vertex scaling for PH metrics on triangulated surfaces.

\begin{definition}[\cite{BPS}]
Suppose $(S,V, \mathcal{T})$ is a connected closed triangulated surface with a PH metric $d: E\rightarrow (0, +\infty)$ and $u: V\rightarrow \mathbb{R}$ is a function defined on the vertices. The vertex scaling of PH metric $d$ by $u$ is defined to be the PH metric $u\ast d$ such that
\begin{equation}\label{discrete conformal equivalent formula}
\sinh\frac{(u\ast d)_{ij}}{2}=\sinh \frac{d_{ij}}{2}e^{u_i+u_j}
\end{equation}
for any edge $\{v_iv_j\}\in E$. The function $u: V\rightarrow \mathbb{R}$ is called a discrete conformal factor.
\end{definition}

As mentioned in \cite{GX2}, the classical combinatorial curvature $K$ is scaling invariant in the Euclidean background geometry
and does not approximate the Gauss curvature on smooth surface as the triangulation of the surface becomes finer and finer.
To overcome these disadvantages, Ge and the first author \cite{GX2} introduced the combinatorial $\alpha$-curvature
for Thurston's Euclidean circle packing metrics.
After that, there are lots of studies on combinatorial $\alpha$-curvature on surfaces and $3$-manifolds \cite{DG,GJ,GX1a,GX2,GX2a,GX1,GX5,X1,X2,Xu1}.
Following \cite{GX2}, we introduce the following combinatorial $\alpha$-curvature for vertex scaling of PH metrics on triangulated surfaces.
\begin{definition}\label{alpha_curvature definition}
Suppose $(S,V, \mathcal{T})$ is a connected closed triangulated surface with a PH metric $d_0$,
$\alpha\in \mathbb{R}$ is a constant and $u: V\rightarrow \mathbb{R}$ is a discrete conformal factor for $d_0$.
The combinatorial $\alpha$-curvature at vertex $v_i\in V$ is defined to be
\begin{equation}\label{alpha_curvature formula}
R_{\alpha, i}=\frac{K_i}{w^{\alpha}_i},
\end{equation}
where $w_i=e^{u_i}$ and $K_i$ is the classical combinatorial curvature at $v_i$ defined by (\ref{K_i}).
\end{definition}
\begin{remark}
$R_{\alpha,i}$ is a parameterized generalization of $K_i$. Specially, if $\alpha=0$, $R_{0,i}$ is the classical combinatorial curvature $K_i$.
Similar to the classical combinatorial curvature $K$, the combinatorial $\alpha$-curvature $R_\alpha$ in Definition \ref{alpha_curvature definition} is independent of the geometric triangulation of the marked
surface $(S,V)$ with a PH metric.
Therefore, the combinatorial $\alpha$-curvature is well-defined for PH metrics discrete conformal to $d_0$.
\end{remark}

One of the remarkable theorems proved by Bobenko-Pinkall-Springborn \cite{BPS}
is that on a triangulated surface with a PH metric, the classical combinatorial curvature determines the discrete conformal factor.
We prove the following global rigidity for the combinatorial $\alpha$-curvature $R_\alpha$ with respect to discrete conformal factors,
which is a parameterized generalization of Bobenko-Pinkall-Springborn's global rigidity.
\begin{theorem}\label{global rigidity theorem}
Suppose $(S,V, \mathcal{T})$ is a connected closed triangulated surface with a PH metric $d_0$,
$\alpha\in \mathbb{R}$ is a constant and $\overline{R}$ is a function defined on the vertices.
If $\alpha\overline{R}\leq 0$,
there exists at most one discrete conformal factor $u$ such that the combinatorial $\alpha$-curvature of $u*d_0$ is $\overline{R}$.
\end{theorem}

An important problem on polyhedral surfaces is about the existence of discrete conformal factors
with prescribed combinatorial curvatures.
Recently, Gu-Luo-Sun-Wu \cite{Gu1},  Gu-Guo-Luo-Sun-Wu \cite{Gu2} and Springborn \cite{Spr}
gave nice answers to this problem. As a corollary, Gu-Luo-Sun-Wu \cite{Gu1} and  Gu-Guo-Luo-Sun-Wu \cite{Gu2} proved the Euclidean
and hyperbolic discrete uniformization theorem for the classical combinatorial curvature on closed marked surfaces respectively.
We prove the following result on the existence of discrete conformal factor for PH metrics with prescribed combinatorial $\alpha$-curvature.

\begin{theorem}\label{discrete uniformization theorem}
Suppose $(S,V)$ is a connected closed marked surface with a PH metric $d_0$, $\alpha$ is a constant and
$\overline{\mathbf{F}}$ is a function defined on the vertices.
Then there exists a PH metric in the discrete conformal class $\mathcal{D}(d_0)$
with combinatorial $\alpha$-curvature $\overline{\mathbf{F}}$ if one of the following three conditions is satisfied
\begin{description}
  \item[(1)] $\alpha>0,\ \chi(S)<0,\ \overline{\mathbf{F}}\leq 0$;
  \item[(2)] $\alpha<0,\ \overline{\mathbf{F}}>0$;
  \item[(3)] (Gu-Guo-Luo-Sun-Wu \cite{Gu2}) $\alpha=0$, $\overline{\mathbf{F}}\in (-\infty, 2\pi)$, $\sum^N_{i=1} \overline{\mathbf{F}}_{i}>2\pi \chi(S)$.
\end{description}
\end{theorem}
Here $\mathcal{D}(d_0)$ is the discrete conformal class introduced by Gu-Guo-Luo-Sun-Wu \cite{Gu2}. See also Definition \ref{Gu. definition}.

As a corollary of Theorem \ref{discrete uniformization theorem},
we have the following result on the existence of PH metrics with constant combinatorial $\alpha$-curvature on connected closed marked surfaces,
which is a parameterized generalization of the hyperbolic discrete uniformization theorem obtained by Gu-Guo-Luo-Sun-Wu for the classical
combinatorial curvature \cite{Gu2}.

\begin{corollary}\label{discrete uniformization theorem constant}
Suppose $(S,V)$ is a connected closed marked surface with a PH metric $d_0$ and $\alpha$ is a constant.
\begin{description}
  \item[(1)] If $\alpha>0$ and $\chi(S)<0$, for any nonpositive constant $\overline{\mathbf{F}}$, there exists a unique PH metric in $\mathcal{D}(d_0)$ with combinatorial $\alpha$-curvature $\overline{\mathbf{F}}$.
  \item[(2)] If $\alpha<0$, for any positive constant $\overline{\mathbf{F}}$, there exists a unique PH metric in $\mathcal{D}(d_0)$ with combinatorial $\alpha$-curvature $\overline{\mathbf{F}}$.
  \item[(3)] (Gu-Guo-Luo-Sun-Wu \cite{Gu2}) If $\alpha=0$, for any constant $\overline{\mathbf{F}}\in (\frac{2\pi\chi(S)}{N}, 2\pi)$,
  there exists a unique PH metric in $\mathcal{D}(d_0)$ with combinatorial $\alpha$-curvature $\overline{\mathbf{F}}$.
\end{description}
\end{corollary}

Another important problem on polyhedral surfaces
is to find effective algorithms searching for polyhedral metrics (PL metrics and PH metrics) with prescribed combinatorial curvatures,
which has lots of applications in practical activities \cite{ZG}.
Chow-Luo's combinatorial Ricci flows for circle packings \cite{Chow-Luo}
and Luo's combinatorial Yamabe flow for vertex scaling \cite{L1}
provide effective algorithms to find polyhedral metrics with prescribed classical combinatorial curvature.
Motivated by \cite{Chow-Luo,GX1,GX5,GX2,L1,Xu1}, we introduce the following combinatorial $\alpha$-Yamabe flow for vertex scaling of PH metrics, which is a parameterized generalization of Luo's combinatorial Yamabe flow in \cite{L1}.

\begin{definition}\label{Ricci flow definition}
Suppose $(S,V, \mathcal{T})$ is a connected closed triangulated surface with a PH metric $d_0$
and $\alpha\in \mathbb{R}$ is a constant. The combinatorial $\alpha$-Yamabe flow is defined to be
\begin{eqnarray}\label{Ricci flow formula}
\begin{cases}
\frac{dw_i}{dt}=-R_{\alpha,i}w_i,\\
w_i(0)=1,
\end{cases}
\end{eqnarray}
where $w_i=e^{u_i}$ and $u: V\rightarrow \mathbb{R}$ is a discrete conformal factor for $d_0$.
\end{definition}

\begin{remark}
When $\alpha=0$, the combinatorial $0$-Yamabe flow (\ref{Ricci flow formula}) is Luo's combinatorial Yamabe flow in \cite{L1}.
One can further generalize the combinatorial $\alpha$-Yamabe flow (\ref{Ricci flow formula}) to the following form to find PH metrics with prescribed combinatorial $\alpha$-curvature $\overline{R}$
\begin{equation}\label{prescribed Ricci flow formula}
\frac{dw_i}{dt}=(\overline{R}_i-R_{\alpha,i})w_i.
\end{equation}
If $\overline{R}_i=0$ for all $i\in V$, (\ref{prescribed Ricci flow formula}) is reduced to (\ref{Ricci flow formula}).
\end{remark}

Combinatorial Calabi flow is another effective combinatorial curvature flow to find polyhedral metrics with prescribed combinatorial curvature. The combinatorial Calabi flow for Thurston's circle packing metrics and for vertex scaling of PL metrics
was introduced by Ge \cite{Ge1,Ge2} and the combinatorial Calabi flow for vertex scaling of PH metrics was introduced by Zhu and the first author \cite{Zhu Xu}.
Motivated by \cite{Ge1,Ge2,Ge3,GX1a,GX1,GX2,Xu1,Zhu Xu}, we introduce the following combinatorial $\alpha$-Calabi flow for vertex scaling of PH metrics.

\begin{definition}\label{Calabi flow definition}
Suppose $(S,V, \mathcal{T})$ is a connected closed triangulated surface with a PH metric $d_0$ and $\alpha\in \mathbb{R}$ is a constant. The combinatorial $\alpha$-Calabi flow is defined to be
\begin{eqnarray}\label{Calabi flow formula}
\begin{cases}
\frac{dw_i}{dt}=\Delta^\mathcal{T}_\alpha R_{\alpha,i}w_i,\\
w_i(0)=1,
\end{cases}
\end{eqnarray}
where $w_i=e^{u_i}$, $u: V\rightarrow \mathbb{R}$ is a discrete conformal factor for $d_0$
and the discrete $\alpha$-Laplace operator $\Delta^\mathcal{T}_\alpha$ is defined by
$$\Delta^\mathcal{T}_\alpha f_i=-\frac{1}{w^\alpha_i}\sum_{j\in V}\frac{\partial K_i}{\partial u_j}f_j$$
for $f: V\rightarrow \mathbb{R}$.
\end{definition}
\begin{remark}
When $\alpha=0$, the combinatorial $0$-Calabi flow (\ref{Calabi flow formula}) is
the combinatorial Calabi flow introduced by Zhu and the first author \cite{Zhu Xu}.
One can also use the following modified combinatorial $\alpha$-Calabi flow to study the prescribed combinatorial $\alpha$-curvature problem
\begin{equation}\label{prescribed Calabi flow formula}
\frac{dw_i}{dt}=\Delta^\mathcal{T}_\alpha(R_\alpha-\overline{R})_iw_i,
\end{equation}
where $\overline{R}$ is a function defined on the vertices.
If $\overline{R}_i=0$ for all $i\in V$, (\ref{prescribed Calabi flow formula}) is reduced to (\ref{Calabi flow formula}).
\end{remark}

Along the combinatorial $\alpha$-flows (combinatorial $\alpha$-Yamabe flow and combinatorial $\alpha$-Calabi flow),
singularities may develop, which correspond to that some triangles degenerate or discrete conformal factors tend to infinity.
To handle the potential singularities along the combinatorial $\alpha$-flows, we do the following surgery on combinatorial $\alpha$-flows by flipping,
the idea of which comes from \cite{L1,Gu1,Gu2}.
Recall that a Delaunay triangulation $\mathcal{T}$ for a PH metric $d$ on a marked surface $(S,V)$ is a geometric triangulation of
$(S,V,d)$ such that for any adjacent triangles $\triangle v_iv_jv_k, \triangle v_iv_jv_l\in F$ sharing an edge $\{v_iv_j\}\in E$,
\begin{equation}\label{Delaunay condition}
\begin{aligned}
\alpha^{ij}_l+\alpha^{ij}_k\leq \alpha^{jk}_i+\alpha^{ik}_j+\alpha^{jl}_i+\alpha^{il}_j.
\end{aligned}
\end{equation}
Please refer to \cite{Gu2,GL} for more information on Delaunay triangulations for PH metrics on marked surfaces.
There is always a Delaunay triangulation $\mathcal{T}$ for a marked surface $(S, V)$ with a PH metric $d$ \cite{Gu2}.
Along the combinatorial $\alpha$-flows on a triangulated surface $(S,V,\mathcal{T})$, if $\mathcal{T}$ is Delaunay in $u(t)\ast d_0$ for $t\in [0,T]$ and not Delaunay in $u(t)\ast d_0$ for $t\in (T,T+\epsilon),\ \epsilon>0$, there exists an edge $\{v_iv_j\}\in E$ such that $\alpha^{jk}_i+\alpha^{ik}_j+\alpha^{jl}_i+\alpha^{il}_j\geq \alpha^{ij}_l+\alpha^{ij}_k$ for $t\in [0,T]$ and $\alpha^{jk}_i+\alpha^{ik}_j+\alpha^{jl}_i+\alpha^{il}_j<\alpha^{ij}_l+\alpha^{ij}_k$ for $t\in (T,T+\epsilon)$. Then we replace the triangulation $\mathcal{T}$ by a new triangulation $\mathcal{T}'$ at time $t=T$ via replacing two triangles $\triangle v_iv_jv_k$ and $\triangle v_iv_jv_l$ adjacent to $\{v_iv_j\}$ by two new triangles $\triangle v_iv_kv_l$ and $\triangle v_jv_kv_l$. This is called a $\mathbf{surgery\  by\  flipping}$ on the triangulation $\mathcal{T}$, which is also an isometry of $(S,V)$ with PH metric $u(T)\ast d_0$.
After the surgery by flipping at time $t=T$, we run the combinatorial $\alpha$-flows on $(S,V, \mathcal{T}')$ with initial metric $u(T)\ast d_0$.
Whenever the Delaunay condition (\ref{Delaunay condition}) is not satisfied along the combinatorial $\alpha$-flows,
we do surgery by flipping along the combinatorial $\alpha$-flows.

\begin{theorem}\label{last theorem}
Suppose $(S,V)$ is a connected closed marked surface with a PH metric $d_0$, $\alpha$ is a constant
and $\overline{\mathbf{F}}$ is a function defined on the vertices.
If one of the conditions $\mathbf{(1)(2)(3)}$ in Theorem \ref{discrete uniformization theorem} is satisfied,
the solution of combinatorial $\alpha$-Yamabe flow and combinatorial $\alpha$-Calabi flow with surgery exist for all time and converge exponentially fast.
\end{theorem}

\begin{remark}
If $\alpha=0$, the longtime existence and convergence of  hyperbolic combinatorial Yamabe flow (0-Yamabe flow) with surgery was proved by Gu-Guo-Luo-Sun-Wu \cite{Gu2} and the longtime existence and convergence of hyperbolic combinatorial Calabi flow (0-Calabi flow) with surgery was proved by Zhu and the first author \cite{Zhu Xu}.
Theorem \ref{last theorem} is a parameterized generalization of these results.
\end{remark}

The paper is organized as follows. In Section \ref{section 2}, we study the combinatorial $\alpha$-curvature on a marked surface with a fixed triangulation and prove Theorem \ref{global rigidity theorem}. We also discuss the stability of combinatorial $\alpha$-flows on triangulated surfaces in Section \ref{section 2}.
In Section \ref{section 3}, we allow the triangulation on the marked surface to be changed by edge flipping and prove Theorem \ref{discrete uniformization theorem} and Theorem \ref{last theorem}.
\\

\textbf{Acknowledgements}\\[8pt]
The authors thank Dr. Tianqi Wu for communications on the details of \cite{Gu2,Gu1}.
The research of the first author is supported by the Fundamental Research Funds for the Central Universities under
grant no. 2042020kf0199.

\section{Combinatorial $\alpha$-curvature and combinatorial $\alpha$-flows on a marked surface with a fixed triangulation}\label{section 2}

\subsection{Rigidity of combinatorial $\alpha$-curvature on triangulated surfaces}
Suppose $(S,V)$ is a marked surface with a fixed triangulation $\mathcal{T}$, $\triangle v_iv_jv_k\in F$ is a triangle and $d_{ij}, d_{ik}, d_{jk}$ defines a PH metric on $\triangle v_iv_jv_k$.
The admissible space $\mathcal{U}^{H,\mathcal{T}}_{ijk}(d)$ of discrete conformal factors for the triangle $\triangle v_iv_jv_k$ with a PH metric $d_{ij}, d_{ik}, d_{jk}$ is defined to be the set of discrete conformal factors $({u_i,u_j,u_k})\in \mathbb{R}^3$ such that the triangle with edge lengths given by (\ref{discrete conformal equivalent formula}) exists in 2-dimensional hyperbolic space $\mathbb{H}^2$, i.e.
\begin{equation}\label{admissible space}
\mathcal{U}^{H,\mathcal{T}}_{ijk}(d)=\{(u_i,u_j,u_k)\in \mathbb{R}^3|(u*d)_{rs}+(u*d)_{rt}>(u*d)_{ts}, \{r,s,t\}=\{i,j,k\} \}.
\end{equation}
Then the admissible space $\mathcal{U}^{H,\mathcal{T}}(d)$ for $(S,V, \mathcal{T})$ with a PH metric $d$ could be written as $\mathcal{U}^{H,\mathcal{T}}(d)=\bigcap_{\triangle v_iv_jv_k\in F}\mathcal{U}^{H,\mathcal{T}}_{ijk}(d)$.
Here and in the following,  we do not distinguish $\mathcal{U}^{H,\mathcal{T}}_{ijk}(d)$
and its natural  preimage in $\mathbb{R}^N$ if it causes no confusion.

For a triangle $\triangle v_iv_jv_k\in F$, we denote $\alpha_i^{jk}$ as the inner angle opposite to the edge $\{v_jv_k\}$.
Based on Leibon's work \cite{GL} on volume of some generalized hyperbolic tetrahedron in $\mathbb{H}^3$,
Bobenko-Pinkall-Springborn \cite{BPS} proved the following result.
\begin{lemma}[\cite{BPS}]\label{constants to be continuous}
Suppose $(S, V, \mathcal{T})$ is a connected closed triangulated surface with a PH metric $d$.
\begin{description}
  \item[(1)] The admissible space $\mathcal{U}^{H,\mathcal{T}}_{ijk}(d)$ is nonempty and simply connected with analytical boundaries in $\mathbb{R}^3$.
  \item[(2)] The matrix $\frac{\partial (\alpha_{i}^{jk}, \alpha_{j}^{ik},\alpha_{k}^{ij})}{\partial (u_i, u_j, u_k)}$ is symmetric and negative definite on $\mathcal{U}^{H,\mathcal{T}}_{ijk}(d)$, which implies the matrix $L$ is symmetric and positive definite on $\mathcal{U}^{H,\mathcal{T}}(d)$, where
\begin{equation}\label{matrix L}
L=(L_{ij})_{N\times N}=\frac{\partial(K_1,...,K_N)}{\partial(u_1,...,u_N)}
 =\left(
 \begin{array}{ccc}
  \frac{\partial K_1}{\partial u_1} & \cdots & \frac{\partial K_1}{\partial u_N} \\
   \vdots  & \ddots & \vdots \\
   \frac{\partial K_N}{\partial u_1} & \cdots & \frac{\partial K_N}{\partial u_N}\\
  \end{array}
 \right).
\end{equation}
  \item[(3)] The inner angles $\alpha_i^{jk},\alpha_j^{ik},\alpha_k^{ij}$ defined for $(u_i,u_j,u_k)\in\mathcal{U}^{H,\mathcal{T}}_{ijk}(d)$
could be extended by constants to be continuous functions $\widetilde{\alpha}_i^{jk},\widetilde{\alpha}_j^{ik},\widetilde{\alpha}_k^{ij}$ defined on $\mathbb{R}^3$.
\end{description}
\end{lemma}

One can also refer to \cite{XZ} for a new proof of Lemma \ref{constants to be continuous}.
We further have the following decomposition for the matrix $L$.
\begin{lemma}\label{L=A+L_B}
Suppose $(S, V, \mathcal{T})$ is a connected closed triangulated surface with a PH metric $d$.
The positive definite matrix $L$ defined by (\ref{matrix L}) could be decomposed to be
$$L=L_A+L_B,$$
where $L_A=\text{diag}\{A_1,...,A_N\}$ is a diagonal matrix with
\begin{equation}\label{A}
A_i=\frac{\partial}{\partial u_i}\left(\sum_{\triangle v_iv_jv_k\in F}\text{Area}(\triangle v_iv_jv_k)\right)
=\sum_{j\sim i}B_{ij}(\cosh d_{ij}-1)
\end{equation}
and $L_B$ is a symmetric matrix defined by
\begin{eqnarray*}
(L_B)_{ij}=
\begin{cases}
\sum_{k\sim i}B_{ik},  &{j=i},\\
-B_{ij}, &{j\sim i},\\
0, &{\text{otherwise}},
\end{cases}
\end{eqnarray*}
with
\begin{equation}\label{B}
B_{ij}
=\frac{\partial \alpha^{jk}_i}{\partial u_j}+\frac{\partial \alpha^{jl}_i}{\partial u_j}
=\frac{1}{\cosh^2 \frac{d_{ij}}{2}}\left(\tan \frac{\alpha^{jk}_i+\alpha^{ik}_j-\alpha^{ij}_k}{2}+\tan \frac{\alpha^{jl}_i+\alpha^{il}_j-\alpha^{ij}_l}{2}\right).
\end{equation}
Here $\triangle v_iv_jv_k$ and $\triangle v_iv_jv_l$ are two adjacent triangles sharing an edge $\{v_iv_j\}$.
Furthermore, if the triangulation $\mathcal{T}$ is Delaunay in the PH metric $d$ on $(S,V)$,
then $A_i>0,\ B_{ij}\geq0$, which implies $L_A$ is positive definite and $L_B$ is positive semi-definite.
\end{lemma}

The proof of Lemma \ref{L=A+L_B} is too long and deviates from the theme of this paper, so we defer it to Appendix \ref{appendix a}.
For Thurston's hyperbolic circle packing metrics on surfaces, Ge and the first author \cite{GX3} obtained a decomposition similar to Lemma \ref{L=A+L_B}.

By Lemma \ref{constants to be continuous},
the Ricci energy function
$$F_{ijk}(u_i,u_j,u_k)=\int^{(u_i,u_j,u_k)}\alpha_i^{jk}du_i+\alpha_j^{ik}du_j+\alpha_k^{ij}du_k$$
for the triangle $\triangle v_iv_jv_k$ is a well-defined locally strictly concave function of $(u_i,u_j,u_k)\in \mathcal{U}^{H,\mathcal{T}}_{ijk}(d)$.
$F_{ijk}(u_i,u_j,u_k)$ was introduced by Bobenko-Pinkall-Springborn \cite{BPS} using the volume function of a generalized
hyperbolic tetrahedron. Here we construct $F_{ijk}(u_i,u_j,u_k)$ directly.
Based on Leibon's work on the volume of a generalized hyperbolic tetrahedron \cite{GT},
Bobenko-Pinkall-Springborn \cite{BPS} further observed that
$F_{ijk}(u_i,u_j,u_k)$ could be extended to be a globally defined concave function on $\mathbb{R}^3$.
In the following, we use Luo's generalization \cite{L3} of Bobenko-Pinkall-Spingborn's extension
to extend the function $F_{ijk}(u_i, u_j, u_k)$.
\begin{definition}[\cite{L3}, Definition 2.3]
A differential 1-form $w=\sum_{i=1}^n a_i(x)dx^i$ in an open subset $U\subset \mathbb{R}^n$ is said to be continuous if each $a_i(x)$ is continuous on $U$. A continuous differential 1-form $w$ is called closed if $\int_{\partial \tau}w=0$ for each triangle $\tau\subset U$.
\end{definition}

\begin{theorem}[\cite{L3}, Corollary 2.6]\label{Luo's convex extention}
Suppose $X\subset \mathbb{R}^n$ is an open convex set and $A\subset X$ is an open subset of $X$ bounded by a real analytic codimension-1 submanifold in $X$. If $w=\sum_{i=1}^na_i(x)dx_i$ is a continuous closed 1-form on $A$ so that $P(x)=\int_a^x w$ is locally convex on $A$ and each $a_i$ can be extended continuous to $X$ by constant functions to a function $\widetilde{a}_i$ on $X$, then $\widetilde{P}(x)=\int_a^x\sum_{i=1}^n\widetilde{a}_i(x)dx_i$ is a $C^1$-smooth
convex function on $X$ extending $P$.
\end{theorem}

Combining Lemma \ref{constants to be continuous} and Theorem \ref{Luo's convex extention}, $F_{ijk}(u_i,u_j,u_k)$ defined on $\mathcal{U}^{H,\mathcal{T}}_{ijk}(d)$ could be extended to be
\begin{equation}\label{extension of Fijk}
\widetilde{F}_{ijk}(u_i,u_j,u_k)=\int^{(u_i,u_j,u_k)}\widetilde{\alpha}_i^{jk}du_i+\widetilde{\alpha}_j^{ik}du_j+\widetilde{\alpha}_k^{ij}du_k,
\end{equation}
which is a $C^1$-smooth concave function defined on $\mathbb{R}^3$.

Using the extension $\widetilde{F}_{ijk}$ of the Ricci energy function $F_{ijk}$,
we can prove the following rigidity for the extended combinatorial $\alpha$-curvature
$\widetilde{R}_\alpha=\frac{\widetilde{K}_i}{w_i^\alpha}=\frac{1}{w_i^\alpha}(2\pi-\sum_{\triangle v_iv_jv_k}\widetilde{\alpha}_i^{jk})$,
which is a generalization of Theorem \ref{global rigidity theorem}.

\begin{theorem}\label{global rigidity theorem context extended}
Suppose $(S,V, \mathcal{T})$ is a connected closed triangulated surface with a PH metric $d_0$,
$\alpha\in \mathbb{R}$ is a constant and $\overline{R}$ is a function defined on the vertices with $\alpha\overline{R}\leq 0$.
If there exist $u_A\in \mathcal{U}^{H,\mathcal{T}}(d_0)$ and $u_B\in \mathbb{R}^N$ with the same extended combinatorial $\alpha$-curvature $\overline{R}$,
then $u_A=u_B$.
\end{theorem}
\proof
Define the following Ricci energy $P(u)$ for $\overline{R}$ by
\begin{equation*}
P(u)=-\sum_{\triangle v_iv_jv_k\in F}F_{ijk}(u_i,u_j,u_k)
+\int^{u}\sum_{i=1}^{N}(2\pi-\overline{R}_iw^{\alpha}_i)du_i.
\end{equation*}
By direct calculations, we have
\begin{equation*}
\nabla_{u_i}P(u)=-\sum_{\triangle v_iv_jv_k\in F}\alpha_i^{jk}+2\pi-\overline{R}_iw^{\alpha}_i
=K_i-\overline{R}_iw^{\alpha}_i, \\
\end{equation*}
and
\begin{equation*}
\mathrm{Hess}_u P=L-\alpha \left(
 \begin{array}{ccc}
  \overline{R}_1w^{\alpha}_1  &  &\\
     & \ddots &  \\
      & & \overline{R}_Nw^{\alpha}_N \\
  \end{array}
 \right).
\end{equation*}
Combining Lemma \ref{constants to be continuous} and the condition $\alpha\overline{R}\leq 0$,
$\mathrm{Hess}_u P$ is positive definite and $P$ is locally strictly convex on $\mathcal{U}^{H,\mathcal{T}}(d_0)$.

By the extension $\widetilde{F}_{ijk}(u_i,u_j,u_k)$ of $F_{ijk}(u_i,u_j,u_k)$ in (\ref{extension of Fijk}),
the Ricci energy function $P(u)$ defined on $\mathcal{U}^{H,\mathcal{T}}(d_0)$ could be extended to be
\begin{equation*}
\widetilde{P}(u)=-\sum_{\triangle v_iv_jv_k\in F}\widetilde{F}_{ijk}(u_i,u_j,u_k)
+\int^{u}\sum_{i=1}^{N}(2\pi-\overline{R}_i\omega^{\alpha}_i)du_i,
\end{equation*}
which is a $C^1$-smooth convex function defined on $\mathbb{R}^N$ and locally strictly convex on $\mathcal{U}^{H,\mathcal{T}}(d_0)\subset \mathbb{R}^N$
with $\nabla_{u_i}\widetilde{P}=\widetilde{K}_i-\overline{R}_iw^{\alpha}_i$.

Set
$$f(t)=\widetilde{P}((1-t)u_A+tu_B),\ t\in [0,1].$$
Then $f(t)$ is a $C^1$ smooth convex function of $t\in [0,1]$ with
$$f'(t)
=\sum_{i=1}^N\nabla_{u_i} \widetilde{P}|_{(1-t)u_A+tu_B}\cdot (u_{B,i}-u_{A,i})
=\sum_{i=1}^N(\widetilde{K}_i-\overline{R}_iw^{\alpha}_i)|_{(1-t)u_A+tu_B}\cdot (u_{B,i}-u_{A,i}).$$
By the assumption that $\widetilde{R}_\alpha(u_A)=\widetilde{R}_\alpha(u_B)=\overline{R}$, we have $f'(0)=f'(1)=0$,
which implies $f'(t)\equiv0$ by the convexity of $f(t)$.
Note that $u_A$ is in the open subset $\mathcal{U}^{H,\mathcal{T}}(d_0)$ of $\mathbb{R}^N$, where $P(u)$ is smooth and
$\mathrm{Hess}_u P$ is positive definite,
there exists $0<\epsilon<1$ such that $(1-t)u_A+tu_B\in \mathcal{U}^{H,\mathcal{T}}(d_0)$
and $f(t)$ is smooth and strictly convex for $t\in[0,\epsilon)$.
By $f'(t)\equiv0$ for $t\in [0, 1]$, we have $f''(t)\equiv0$ for $t\in [0, \epsilon)$.
Combining with
$$f''(t)=(u_B-u_A)^T\cdot \mathrm{Hess}_u P|_{(1-t)u_A+tu_B} \cdot (u_B-u_A)$$
for $t\in [0, \epsilon)$ and $\mathrm{Hess}_u P$ is positive definite for $u\in \mathcal{U}^{H,\mathcal{T}}(d_0)$,
this implies $u_A=u_B$.
\qed

\subsection{Combinatorial $\alpha$-flows on triangulated surfaces}
The combinatorial $\alpha$-Yamabe flow (\ref{prescribed Ricci flow formula})
and the combinatorial $\alpha$-Calabi flow (\ref{prescribed Calabi flow formula}) are ODE systems
with smooth coefficients.
Therefore, the solutions always exist locally around the initial time $t=0$.
We further have the following result on the longtime existence and convergence for the solutions of combinatorial $\alpha$-flows.

\begin{theorem}\label{alpha flows with small initial energy}
Suppose $(S,V,\mathcal{T})$ is a connected closed triangulated surface with a PH metric $d_0$, $\alpha\in \mathbb{R}$ is a constant and
$\overline{R}$ is a function defined on $V$.
If the solutions of combinatorial $\alpha$-flows converge, there exist PH metrics with combinatorial $\alpha$-curvature $\overline{R}$.
Furthermore, if there is a PH metric $d^{\ast}=u^{\ast}\ast d_0$ with combinatorial $\alpha$-curvature $\overline{R}$ and $\alpha \overline{R}\leq 0$, there exists a constant $\delta >0$ such that if $||R_\alpha(u(0))-R_\alpha(u^\ast)||<\delta$, then the solutions of combinatorial $\alpha$-flows exist for all time and converge exponentially fast to $u^\ast$.
\end{theorem}
\proof
Suppose $w(t)$ is a solution of the combinatorial $\alpha$-Yamabe flow (\ref{prescribed Ricci flow formula}),
then $u(t)=\ln w(t)$ is a solution of the flow
$\frac{du_i}{dt}=\overline{R}_i-R_{\alpha,i}.$
If $w^\ast:=w(+\infty)=\lim_{t\rightarrow +\infty}w(t)$ exists,
then $u^\ast:=u(+\infty)=\ln w(+\infty)$ exists in $\mathcal{U}^{H,\mathcal{T}}(d_0)$ and $R_\alpha(u^\ast)=\lim_{t\rightarrow +\infty}R_\alpha(u(t))$ exists.
Furthermore, there exists a sequence $\xi_n\in(n,n+1)$ such that
$$u_i(n+1)-u_i(n)=u'_i(\xi_n)=\overline{R}_i-R_{\alpha,i}(u(\xi_n))\rightarrow 0\  \text{as} \ n\rightarrow +\infty,$$
which implies $R_{\alpha,i}(u^\ast)=\overline{R}_i$ for all $i\in V$
and $u^\ast\ast d_0$ is a PH metric with combinatorial $\alpha$-curvature $\overline{R}$.
Similarly, if the combinatorial $\alpha$-Calabi flow (\ref{prescribed Calabi flow formula}) converges, there exists a PH metric with combinatorial $\alpha$-curvature $\overline{R}$.

Suppose there exists a PH metric $d^{\ast}=u^{\ast}\ast d_0$ with combinatorial $\alpha$-curvature $\overline{R}$.
For the combinatorial $\alpha$-Yamabe flow (\ref{prescribed Ricci flow formula}),
set $\Gamma(u)=\overline{R}-R_{\alpha}$. By direct calculations, we have
$$D\Gamma|_{u=u^*}=-\Sigma^{-\alpha}L+\alpha \Lambda =-\Sigma^{-\frac{\alpha}{2}}
(\Sigma^{-\frac{\alpha}{2}}L\Sigma^{-\frac{\alpha}{2}}-\alpha\Lambda)\Sigma^{\frac{\alpha}{2}},$$
where $\Sigma=\text{diag}\{w_1,w_2,...,w_N\}$, $\Lambda=\text{diag}\{\overline{R}_1,\overline{R}_2,...,\overline{R}_N\}$.
By the condition $\alpha \overline{R}\leq 0$, $D\Gamma|_{u=u^*}$ has $N$ negative eigenvalues,
which implies $u^*$ is a local attractor of (\ref{prescribed Ricci flow formula}).
Then the conclusion follows from Lyapunov Stability Theorem (\cite{Pontryagin}, Chapter 5).

Similarly, for the combinatorial $\alpha$-Calabi flow (\ref{prescribed Calabi flow formula}), set $\Gamma(u)=\Delta^\mathcal{T}_\alpha(R_\alpha-\overline{R})$. By direct calculations, we have 
\begin{equation*}
\begin{aligned}
D\Gamma|_{u=u^*}
&=-\Sigma^{-\alpha}L\Sigma^{-\alpha}L+\alpha\Sigma^{-\alpha}L\Lambda \\
&=-\Sigma^{-\frac{\alpha}{2}}\left(\Sigma^{-\frac{\alpha}{2}}L\Sigma^{-\alpha}L\Sigma^{-\frac{\alpha}{2}}-\alpha
\Sigma^{-\frac{\alpha}{2}} L\Sigma^{-\frac{\alpha}{2}} \Lambda \right)\Sigma^{\frac{\alpha}{2}}\\
&=-\Sigma^{-\frac{\alpha}{2}}\left(Q^2-\alpha Q \Lambda\right)\Sigma^{\frac{\alpha}{2}}\\
&=-\Sigma^{-\frac{\alpha}{2}}Q^{\frac{1}{2}}\left(Q^2+Q^{\frac{1}{2}}(-\alpha \Lambda)Q^{\frac{1}{2}}\right)Q^{-\frac{1}{2}}\Sigma^{\frac{\alpha}{2}},
\end{aligned}
\end{equation*}
where $Q=\Sigma^{-\frac{\alpha}{2}} L\Sigma^{-\frac{\alpha}{2}}$ is a symmetric and positive definite matrix.
By the condition $\alpha \overline{R}\leq 0$, $D\Gamma|_{u=u^*}$ has $N$ negative eigenvalues, which implies $u^\ast$ is a local attractor of (\ref{prescribed Ricci flow formula}). Then the conclusion follows from Lyapunov Stability Theorem (\cite{Pontryagin}, Chapter 5).
\qed

Theorem \ref{alpha flows with small initial energy} gives the longtime existence and convergence for
the solutions of combinatorial $\alpha$-flows for initial PH metrics with small initial energy.
But for general initial PH metrics, the combinatorial $\alpha$-flows may develop singularities, including that some triangles degenerate and the conformal factors tend to infinity along the corresponding combinatorial $\alpha$-flows.
To handle the possible singularities along the combinatorial $\alpha$-flows, we can do surgery by flipping on the $\alpha$-flows.

\section{Combinatorial $\alpha$-curvature and combinatorial $\alpha$-flows on discrete Riemannian surfaces}\label{section 3}
\subsection{Gu-Guo-Luo-Sun-Wu's work on discrete uniformization theorem}
To analyze the behavior of the combinatorial $\alpha$-flows with surgery,
we need the discrete conformal theory for PH metrics established by Gu-Guo-Luo-Sun-Wu \cite{Gu2}.
In this subsection, we briefly recall the theory. Please refer to \cite{Gu2} for more details.
\begin{definition}(\cite{Gu2}, Definition 1)\label{Gu. definition}
Two PH metrics $d,d'$ on a connected closed marked surface $(S,V)$ are discrete conformal
if there exist a sequence of PH metrics $d_1=d,d_2,...,d_m=d'$ on $(S,V)$ and triangulations $\mathcal{T}_1,...,\mathcal{T}_m$ of $(S,V)$ satisfying
\begin{description}
  \item[(a)] each $\mathcal{T}_i$ is Delaunay in $d_i$,
  \item[(b)] if $\mathcal{T}_i=\mathcal{T}_{i+1}$, there exists a function $u: V\rightarrow \mathbb{R}$, called a conformal factor, so that if $e$ is an edge in $\mathcal{T}_i$ with end points $v$ and $v'$, then the lengths $x_{d_i}(e)$ and $x_{d_{i+1}}(e)$ of $e$ in metrics $d_i$ and $d_{i+1}$ are related by $$\sinh\frac{x_{d_{i+1}}(e)}{2}=\sinh \frac{x_{d_i}(e)}{2}e^{u(v)+u(v')},$$
  \item[(c)] if $\mathcal{T}_i\neq\mathcal{T}_{i+1}$, then $(S,d_i)$ is isometric to $(S,d_{i+1})$ by an isometry homotopic to the identity in $(S,V)$.
\end{description}
\end{definition}
The discrete conformal class of a PH metric $d$ on $(S,V)$ is called a discrete Riemannian surface, which is denoted by $\mathcal{D}(d)$ \cite{Gu2}.

\begin{lemma}(\cite{Gu2}, Proposition 16)\label{Gu. diagonal switch}
If $\mathcal{T}$ and $\mathcal{T}'$ are Delaunay triangulations of a PH metric $d$ on a closed marked surface $(S,V)$, there exist a sequence of Delaunay triangulations $\mathcal{T}_1=\mathcal{T},\mathcal{T}_2,..., \mathcal{T}_k=\mathcal{T}'$ of $d$ so that $\mathcal{T}_{i+1}$ is obtained from $\mathcal{T}_i$ by a diagonal switch.
\end{lemma}
The diagonal switch in Lemma \ref{Gu. diagonal switch} is the surgery by flipping described in the introduction.
Using the new definition of discrete conformality, Gu-Guo-Luo-Sun-Wu \cite{Gu2} proved the following
discrete conformal uniformization theorem for vertex scaling of PH metrics.

\begin{theorem}(\cite{Gu2}, Theorem 3)\label{Gu. discrete uniformization theorem}
Suppose $(S,V)$ is a connected closed marked surface and $d$ is a PH metric on $(S,V)$.
Then for any $K^\ast : V\rightarrow (-\infty, 2\pi)$ with $\sum_{v\in V}K^\ast(v)>2\pi \chi(S)$, there exists a unique PH metric $d'$ on $(S,V)$ so that $d'$ is discrete conformal to $d$ and the discrete curvature of $d'$ is $K^\ast$. Furthermore, the discrete Yamabe flow with surgery associated to curvature $K^\ast$ having initial value $d$ converges to $d'$ linearly fast.
\end{theorem}
Denote the Teichim\"{u}ller space of all PH metrics on $(S,V)$ by $T_{hp}(S,V)$ and decorated Teichim\"{u}ller space of all equivalence class of decorated hyperbolic metrics on $S-V$ by $T_D(S-V)$. To prove Theorem \ref{Gu. discrete uniformization theorem}, Gu-Guo-Luo-Sun-Wu
 proved the following result.

\begin{theorem}(\cite{Gu2}, Theorem 22, Corollary 24)\label{Gu.  Theorem 22}
There is a $C^1$-diffeomorphism $\mathbf{A}: T_{hp}(S,V)\rightarrow T_D(S,V)$ between $T_{hp}(S,V)$ and $T_D(S,V)$. Furthermore, the space $\mathcal{D}(d)\subset T_{hp}(S,V)$ of all equivalence classes of PH metrics discrete conformal to $d$ is $C^1$-diffeomorphism to $\{p\}\times \mathbb{R}^{V}_{>0}$ under the diffeomorphism $\mathbf{A}$, where $p$ is the unique hyperbolic metric on $S-V$ determined by the PH metric $d$ on $(S,V)$.
\end{theorem}

Theorem \ref{Gu.  Theorem 22} implies that the union of the admissible spaces $\mathcal{U}^{H,\mathcal{T}}_D(d')$ of conformal factors such that $\mathcal{T}$ is Delaunay for $d'\in \mathcal{D}(d)$ is $\mathbb{R}^n$.

Set $u_i=\ln w_i,\ 1\leq i\leq n$, for $w=(w_1,w_2,...,w_n)\in \mathbb{R}^n_{>0}$ and $\mathbf{P}=\{x\in(-\infty,2\pi)^n|\sum^n_{i=1}x_i>2\pi\chi(S)\}$. Using the map $\mathbf{A}$, Gu-Guo-Luo-Sun-Wu \cite{Gu2} defined the curvature map
\begin{equation}\label{F definition}
\begin{aligned}
\mathbf{F} :\  &\mathbb{R}^n\rightarrow \mathbf{P}\\
&u\mapsto K_{\mathbf{A}^{-1}(p,w(-2u))}
\end{aligned}
\end{equation}
and proved the following property of $\mathbf{F}$.

\begin{proposition}(\cite{Gu2})\label{W}
There exists a $C^2$-smooth strictly convex function $W:\mathbb{R}^n\rightarrow \mathbb{R}$
defined by $W(u)=\int^u_{0}\sum_{i=1}^N\mathbf{F}_idu_i$ so that its
gradient $\nabla W$ is $\mathbf{F}$.
\end{proposition}

\begin{remark}\label{PH metric determine conformal factor}
Proposition \ref{W} implies that $\mathbf{F}$ defined on $\mathbb{R}^n$ is a $C^1$-extension of the combinatorial
curvature $K$ defined on the space of conformal factors $\mathcal{U}^{H,\mathcal{T}}_D(d')$ for $d'\in \mathcal{D}(d)$.
Furthermore, the discrete conformal factor $u$ is uniquely determined by the curvature map $\mathbf{F}$ and then determined
by the PH metric in the discrete conformal class $\mathcal{D}(d)$ \cite{Gu2}.
\end{remark}

Note that the discrete $\alpha$-Laplace operator $\Delta^\mathcal{T}_\alpha$ is independent of the Delaunay triangulations of a PH metric, then the discrete $\alpha$-Laplace operator $\Delta^\mathcal{T}_\alpha$ could be extend to the following operator $\Delta_\alpha$ defined on $\mathbb{R}^n$, which is the space of discrete conformal factors for the discrete conformal class $\mathcal{D}(d)$.

\begin{definition}\label{Delta alpha}
Suppose $(S, V)$ is a connected closed marked surface with a PH metric $d_0$.
For a function $f : V \rightarrow \mathbb{R}$ on the vertices, the hyperbolic discrete $\alpha$-Laplace operator of $d\in \mathcal{D}(d_0)$ on $(S, V)$ is defined to be the map
\begin{equation*}
\begin{aligned}
\Delta_\alpha :\  &\mathbb{R}^{V}\rightarrow \mathbb{R}^{V}\\
&f\mapsto \Delta_\alpha f,
\end{aligned}
\end{equation*}
where the value of $\Delta_\alpha f$ at vertex $v_i$ is
\begin{equation}\label{alpha laplacian definition}
\Delta_\alpha f_i=-\frac{1}{w^\alpha_i}\sum_{j\in V}\frac{\partial \mathbf{F}_i}{\partial u_j}f_j=-\frac{1}{w^\alpha_i}(\widetilde{L}f)_i
\end{equation}
with $\widetilde{L}_{ij}=\frac{\partial \mathbf{F}_i}{\partial u_j}$ being an extension of $L_{ij}=\frac{\partial K_i}{\partial u_j}$.
\end{definition}

\begin{remark}\label{remark laplace decomp Delaunay}
By Remark \ref{PH metric determine conformal factor}, the discrete conformal factor $u\in \mathbb{R}^N$ is uniquely
determined by the PH metric $d\in \mathcal{D}(d_0)$ and the discrete $\alpha$-Laplace operator in Definition \ref{Delta alpha} is well-defined.
By Proposition \ref{W}, $\Delta_\alpha$ is continuous and piecewise smooth on $\mathbb{R}^n$ as a matrix-valued function of $u$. By Lemma \ref{L=A+L_B}, we have
\begin{equation}\label{Delta alpha f_i}
\Delta_\alpha f_i=\sum_{j\sim i}\frac{B_{ij}}{w^\alpha_i}(f_j-f_i)-\frac{A_i}{w^\alpha_i}f_i
\end{equation}
with $A_i>0$ and $B_{ij}\geq0$.
\end{remark}

\subsection{Rigidity of combinatorial $\alpha$-curvature on discrete Riemannian surfaces}
Motivated by the relationship of $\mathbf{F}$ defined on a discrete Riemannian surface and the combinatorial
curvature $K$ defined on a triangulated surface,
we introduce the following definition of
combinatorial $\alpha$-curvature $\mathbf{F}_{\alpha}$ for PH metrics on a discrete Riemannian surface,
which is an extension of the combinatorial $\alpha$-curvature $R_\alpha$ defined on a triangulated surface.
\begin{definition}\label{F alpha_curvature definition}
Suppose $(S,V)$ is a connected closed marked surface with a PH metric $d_0$, $\alpha \in \mathbb{R}$ is a constant and
$\mathbf{F}$ is the curvature map defined by (\ref{F definition}). The combinatorial $\alpha$-curvature for $d\in\mathcal{D}(d_0)$ is defined to be
\begin{equation}\label{F alpha formula}
\mathbf{F}_{\alpha, i}=\frac{\mathbf{F}_i}{w^{\alpha}_i},
\end{equation}
where $w_i=e^{u_i}$ and $u: V\rightarrow \mathbb{R}$ is the unique discrete conformal factor determined by
the PH metric $d\in \mathcal{D}(d_0)$.
\end{definition}
By Remark \ref{PH metric determine conformal factor},
the combinatorial $\alpha$-curvature $\mathbf{F}_\alpha$ in Definition \ref{F alpha_curvature definition} is well-defined for $d\in \mathcal{D}(d_0)$.
Denote the space of conformal factors by $\mathcal{U}(d)$ for a PH metric $d$ on the marked surface $(S,V)$,
which is $\mathbb{R}^N$ by Theorem \ref{Gu.  Theorem 22}.
Analogous to Theorem \ref{global rigidity theorem}, we have the following global rigidity for the combinatorial $\alpha$-curvature $\mathbf{F}_\alpha$
with respect to discrete conformal factors on discrete Riemannian surfaces.

\begin{theorem}\label{F global rigidity theorem}
Suppose $(S,V)$ is a closed marked surface with a PH metric $d_0$, $\alpha\in \mathbb{R}$ is a constant and
$\overline{F}$ is a function defined on $V$.
If $\alpha\overline{F}\leq 0$, there exists at most one discrete conformal factor $u^*\in \mathcal{U}(d_0)$ such that $\mathbf{A}^{-1}(p,w(-2u^\ast))\in \mathcal{D}(d_0)$ has the combinatorial $\alpha$-curvature $\overline{F}$.
\end{theorem}
\proof
By Proposition \ref{W}, define the following Ricci energy function
\begin{equation*}
W_\alpha(u)=W(u)-\int_{0}^{u}\sum_{i=1}^{N}\overline{F}_iw^{\alpha}_idu_i
=\int_0^u\sum_{i=1}^N\mathbf{F}_idu_i-\int_{0}^{u}\sum_{i=1}^{N}\overline{F}_iw^{\alpha}_idu_i.
\end{equation*}
Then $W_\alpha(u)$ is a $C^2$-smooth function defined on $\mathbb{R}^n$ with
$\nabla_{u_i}W_\alpha=\nabla_{u_i}W-\overline{F}_iw^\alpha_i=\mathbf{F}_i-\overline{F}_iw^\alpha_i.$
Therefore, for $u^*\in \mathcal{U}(d)$, $\mathbf{A}^{-1}(p,w(-2u^\ast))\in \mathcal{D}(d)$ has the combinatorial $\alpha$-curvature $\overline{F}$ if and only if $\nabla_{u_i}W_\alpha(u^\ast)=0$. By direct calculations, we have
\begin{equation*}
\mathrm{Hess} W_\alpha=\widetilde{L}-\alpha \left(
 \begin{array}{ccc}
  \overline{F}_1w^{\alpha}_1  &  &\\
     & \ddots &  \\
      & & \overline{F}_Nw^{\alpha}_N \\
  \end{array}
 \right).
\end{equation*}
Due to $\alpha\overline{F}\leq 0$, $\mathrm{Hess} W_\alpha$ is positive definite and $W_\alpha$ is strictly convex on
$\mathbb{R}^n$ by Proposition \ref{W}. Then the rigidity follows from the following well-known result from analysis.
\begin{lemma}\label{1}
If $W : \Omega \rightarrow \mathbb{R}$ is a $C^1$-smooth  strictly convex function on an open convex set $\Omega \subset \mathbb{R}^m$, then its gradient $\nabla W : \Omega \rightarrow \mathbb{R}^m$ is an embedding.
\end{lemma}
\qed

\subsection{Combinatorial $\alpha$-Yamabe flow with surgery}
Using the discrete conformal theory established by Gu-Guo-Luo-Sun-Wu \cite{Gu2},
the combinatorial $\alpha$-Yamabe flow with surgery can be defined as follows.
\begin{definition}\label{Ricci flow with surgery definition}
Suppose $(S,V)$ is a connected closed marked surface with a PH metric $d_0$ and $\alpha\in \mathbb{R}$ is a constant. The combinatorial $\alpha$-Yamabe flow with surgery is defined to be
\begin{eqnarray}\label{Ricci flow with surgery formula}
\begin{cases}
\frac{du_i}{dt}=-\mathbf{F}_{\alpha,i},\\
u_i(0)=0.
\end{cases}
\end{eqnarray}
\end{definition}
One can modify the combinatorial $\alpha$-Yamabe flow with surgery (\ref{Ricci flow with surgery formula}) to the following form
\begin{equation}\label{prescribed Ricci flow with surgery formula}
\frac{du_i}{dt}=\overline{\mathbf{F}}_i-\mathbf{F}_{\alpha,i}
\end{equation}
to study prescribed curvature problems, where $\overline{\mathbf{F}}$ is a function defined on $V$.
\begin{theorem}\label{equivalent theorem 1}
Suppose $(S,V)$ is a connected closed marked surface with a PH metric $d_0$, $\alpha\in \mathbb{R}$ is a constant
and $\overline{\mathbf{F}}: V\rightarrow \mathbb{R}$ is a function defined on $V$ with $\alpha \overline{\mathbf{F}}\leq 0$.
Then there exists a PH metric in the discrete conformal class $\mathcal{D}(d_0)$
with combinatorial $\alpha$-curvature $\overline{\mathbf{F}}$ if and only if
the solution of modified combinatorial $\alpha$-Yamabe flow with surgery (\ref{prescribed Ricci flow with surgery formula}) exists for all time and converges exponentially fast to a PH metric with combinatorial $\alpha$-curvature $\overline{\mathbf{F}}$.
\end{theorem}
\proof
Suppose the solution $u(t)$ of the modified combinatorial $\alpha$-Yamabe flow with surgery (\ref{prescribed Ricci flow with surgery formula}) converges to $u^{\ast}\in \mathcal{U}(d_0)$ as $t\rightarrow +\infty$, then $\mathbf{F}_\alpha(u^*)=\lim_{t\rightarrow +\infty}\mathbf{F}_\alpha(u(t))$ by the $C^1$-smoothness of $\mathbf{F}_\alpha$.
Furthermore, there exists a sequence $\xi_n\in(n,n+1)$ such that
\begin{equation*}
u_i(n+1)-u_i(n)=u'_i(\xi_n)=\overline{\mathbf{F}}_i-\mathbf{F}_{\alpha,i}(u(\xi_n))\rightarrow 0,\ \text{as}\  n\rightarrow +\infty,
\end{equation*}
which implies $\mathbf{F}_{\alpha,i}(u^*)=\lim_{n\rightarrow +\infty}\mathbf{F}_{\alpha,i}(u(\xi_n))=\overline{\mathbf{F}}_i$ for all $v_i\in V$
and $u^\ast\in \mathcal{U}(d_0)$ is a discrete conformal factor with combinatorial $\alpha$-curvature $\overline{\mathbf{F}}$.

Conversely, suppose there exists a PH metric $u^\ast$ in the discrete conformal class $\mathcal{D}(d_0)$ with combinatorial $\alpha$-curvature $\overline{\mathbf{F}}$. Set
\begin{equation*}
W_\alpha(u)=W(u)-\int_{u^*}^{u}\sum_{i=1}^{N}\overline{\mathbf{F}}_iw^{\alpha}_idu_i
=\int^u_{u^*}\sum_{i=1}^N\mathbf{F}_idu_i-\int_{u^*}^{u}\sum_{i=1}^{N}\overline{\mathbf{F}}_iw^{\alpha}_idu_i.
\end{equation*}
Then $W_\alpha(u)$ is a $C^2$-smooth strictly convex function defined on $\mathbb{R}^n$ by Proposition \ref{W} and the condition $\alpha \overline{\mathbf{F}}\leq 0$.
Furthermore, $W_\alpha(u^*)=0,\ \nabla W_\alpha(u^*)=0$,
which implies $W_\alpha(u)\geq W_\alpha(u^\ast)=0$ and $\lim_{u\rightarrow \infty}W_\alpha(u)=+\infty$ by the convexity of $W_\alpha(u)$.
On can also refer to Lemma 4.6 in \cite{GX5} for a proof for this fact.
By direct calculations, we have
\begin{equation*}
\frac{dW_\alpha(u(t))}{dt}
=\sum^N_{i=1}\frac{\partial W_\alpha}{\partial u_i}\frac{du_i}{dt}
=\sum^N_{i=1}(\mathbf{F}_{\alpha,i}-\overline{\mathbf{F}}_i)w^\alpha_i
(\overline{\mathbf{F}}_i-\mathbf{F}_{\alpha,i})
=-\sum^N_{i=1}(\overline{\mathbf{F}}_i-\mathbf{F}_{\alpha,i})^2w^\alpha_i\leq 0,
\end{equation*}
which implies $0\leq W_\alpha(u(t))\leq W_\alpha(u(0))$.
Therefore, $\{u(t)\}\subset\subset \mathbb{R}^N$ by $\lim_{u\rightarrow \infty}W_\alpha(u)=+\infty$,
which implies the solution of modified combinatorial $\alpha$-Yamabe flow with surgery (\ref{prescribed Ricci flow with surgery formula}) exists for all time and $W_\alpha(u(t))$ converges.
Moreover, there exists a sequence $\xi_n\in(n,n+1)$ such that as $n\rightarrow +\infty$,
\begin{equation*}
W_\alpha(u(n+1))-W_\alpha(u(n))=(W_\alpha(u(t))'|_{\xi_n}=\nabla W_\alpha\cdot\frac{du_i}{dt}|_{\xi_n}
=-\sum^N_{i=1}(\overline{\mathbf{F}}_{i}-\mathbf{F}_{\alpha,i})^2w^\alpha_i|_{\xi_n} \rightarrow 0.
\end{equation*}
Then $\lim_{n\rightarrow +\infty}\mathbf{F}_{\alpha,i}(u(\xi_n))=\overline{\mathbf{F}}_i=\mathbf{F}_{\alpha,i}(u^*)$ for all $v_i\in V$.
By $\{u(t)\}\subset\subset \mathbb{R}^N$, there exists $\overline{u}\in \mathbb{R}^N$ and
a subsequence of $\{u(\xi_n)\}$, still denoted as $\{u(\xi_n)\}$ for simplicity,
such that $\lim_{n\rightarrow \infty}u(\xi_n)=\overline{u}$, which implies
$\mathbf{F}_{\alpha,i}(\overline{u})=\lim_{n\rightarrow +\infty}\mathbf{F}_{\alpha,i}(u(\xi_n))=\mathbf{F}_{\alpha,i}(u^*)$.
This further implies $\overline{u}=u^*$ by Theorem \ref{F global rigidity theorem}.
Therefore, $\lim_{n\rightarrow \infty}u(\xi_n)=u^*$.

Set $\Gamma(u)=\overline{\mathbf{F}}-\mathbf{F}_{\alpha}$, then $D\Gamma|_{u=u^\ast}$ has $N$ negative eigenvalues, which implies that $u^\ast$ is a local attractor of (\ref{prescribed Ricci flow with surgery formula}). Then the conclusion follows from Lyapunov Stability Theorem (\cite{Pontryagin}, Chapter 5).
\qed

\subsection{Existence of PH metrics with the prescribed $\alpha$-curvature}
Set $M_{\alpha,i}=\mathbf{F}_{\alpha,i}-\overline{\mathbf{F}}_i$.
Then the modified combinatorial Yamabe flow with surgery (\ref{prescribed Ricci flow with surgery formula}) could be written as $$\frac{du_i}{dt}=-M_{\alpha,i}.$$
We have the following evolution equation for $M_{\alpha,i}$ along the modified combinatorial $\alpha$-Yamabe flow with surgery (\ref{prescribed Ricci flow with surgery formula}).
\begin{lemma}
Along the modified combinatorial $\alpha$-Yamabe flow with surgery (\ref{prescribed Ricci flow with surgery formula}), $M_{\alpha,i}$ evolves according to
\begin{equation}\label{the evolution of M alpha i}
\frac{dM_{\alpha,i}}{dt}
=\sum_{j\sim i}\frac{B_{ij}}{w^\alpha_i}(M_{\alpha,j}-M_{\alpha,i}) -\frac{A_i}{w^\alpha_i}M_{\alpha,i}+\alpha M_{\alpha,i}(M_{\alpha,i}+\overline{\mathbf{F}}_i),
\end{equation}
where $A_i>0$, $B_{ij}\geq0.$
\end{lemma}
\proof
By the chain rules,
\begin{equation*}
\begin{aligned}
\frac{dM_{\alpha,i}}{dt}
=&\sum^N_{j=1}\frac{\partial \mathbf{F}_{\alpha,i}}{\partial u_j}\frac{du_j}{dt}\\
=&\sum^N_{j=1}(\frac{1}{w^\alpha_i}\frac{\partial \mathbf{F}_i}{\partial u_j}-\frac{\mathbf{F}_i}{w^\alpha_i}\alpha \delta_{ij})\cdot (-M_{\alpha,j})\\
=&(\Delta_\alpha M_\alpha)_i+\alpha M_{\alpha,i}(M_{\alpha,i}+\overline{\mathbf{F}}_i)\\
=&\sum_{j\sim i}\frac{B_{ij}}{w^\alpha_i}(M_{\alpha,j}-M_{\alpha,i}) -\frac{A_i}{w^\alpha_i}M_{\alpha,i}+\alpha M_{\alpha,i}(M_{\alpha,i}+\overline{\mathbf{F}}_i),
\end{aligned}
\end{equation*}
where (\ref{Delta alpha f_i}) is used in the last line.
$A_i>0$ and $B_{ij}\geq0$ follows from Lemma \ref{L=A+L_B}.
\qed

To analyze the longtime behavior of the combinatorial $\alpha$-Yamabe flow with surgery, we need the following
discrete maximum principle obtained in \cite{GX2}.

\begin{theorem}(Maximum Principle)\label{maximum principle}
Let $f:V\times [0,T)\rightarrow \mathbb{R}$ be a $C^1$ function such that
\begin{equation*}
\frac{\partial f_i}{\partial t}\geq \sum_{j\sim i}a_{ij}(f_j-f_i)+\Phi_i(f_i),\ \forall(v_i,t)\in V\times [0,T),
\end{equation*}
where $a_{ij}\geq 0$ and $\Phi_i: \mathbb{R}\rightarrow \mathbb{R}$ is a local Lipschitz function. Suppose there exists $C_1\in \mathbb{R}$ such that $f_i(0)\geq C_1$ for all $v_i\in V$. Let $\varphi$ be the solution to the associated ODE
\begin{eqnarray*}
\begin{cases}
\frac{d\varphi}{dt}=\Phi_i(\varphi)\\
\varphi(0)=C_1,
\end{cases}
\end{eqnarray*}
then $f_i(t)\geq \varphi(t)$ for all $(v_i,t)\in V\times [0,T)$ such that $\varphi(t)$ exists.

Similarly, suppose $f:V\times [0,T)\rightarrow \mathbb{R}$ be a $C^1$ function such that
\begin{equation*}
\frac{\partial f_i}{\partial t}\leq \sum_{j\sim i}a_{ij}(f_j-f_i)+\Phi_i(f_i),\ \forall(v_i,t)\in V\times [0,T).
\end{equation*}
Suppose there exists $C_2\in \mathbb{R}$ such that $f_i(0)\leq C_2$ for all $v_i\in V$. Let $\psi$ be the solution to the associated ODE
\begin{eqnarray*}
\begin{cases}
\frac{d\psi}{dt}=\Phi_i(\psi)\\
\psi(0)=C_2,
\end{cases}
\end{eqnarray*}
then $f_i(t)\leq \psi(t)$ for all $(v_i,t)\in V\times [0,T)$ such that $\psi(t)$ exists.
\end{theorem}

\begin{remark}
In Theorem \ref{maximum principle}, $\sum_{j\sim i}a_{ij}(f_j-f_i)$ is actually the classical discrete Laplace operator $\Delta f_i$
used in \cite{Chung}. In order to avoid confusion with the $\alpha$-Laplace operator $\Delta_\alpha$ in Definition \ref{Delta alpha},
here we use the notation $\sum_{j\sim i}a_{ij}(f_j-f_i)$ rather than $\Delta f_i$.
\end{remark}
As a direct application of Theorem \ref{maximum principle} to (\ref{the evolution of M alpha i}), we have the following property for $M_{\alpha}$.
\begin{corollary}\label{application of maximum principle}
If $M_{\alpha,i}(0)\leq 0$ for all $v_i\in V$, then $M_{\alpha,i}(t)\leq 0$ for all $v_i\in V$ along the modified combinatorial $\alpha$-Yamabe flow (\ref{prescribed Ricci flow with surgery formula}). If $M_{\alpha,i}(0)\geq 0$ for all $v_i\in V$, then $M_{\alpha,i}(t)\geq 0$ for all $v_i\in V$ along the modified combinatorial $\alpha$-Yamabe flow (\ref{prescribed Ricci flow with surgery formula}).
\end{corollary}

As another application of Theorem \ref{maximum principle}, we have the following existence of PH metrics with
prescribed combinatorial $\alpha$-curvature.

\begin{lemma}\label{Lemma 1}
Suppose $(S,V)$ is a connected closed marked surface with a PH metric $d_0$, $\alpha>0$ is a constant and $\overline{\mathbf{F}}: V\rightarrow \mathbb{R}$ is a negative constant function defined on $V$.
If there exists $u^\ast\in \mathcal{U}(d_0)$ such that $\overline{\mathbf{F}}<\mathbf{F}_{\alpha,i}(u^\ast)<0$ for all $v_i\in V$,
then there exists a PH metric in the discrete conformal class $\mathcal{D}(d_0)$ with combinatorial $\alpha$-curvature $\overline{\mathbf{F}}$.
\end{lemma}
\proof
Taking $u^*$ as the initial value of
the modified combinatorial $\alpha$-Yamabe flow with surgery (\ref{prescribed Ricci flow with surgery formula}).
We shall prove that the solution of (\ref{prescribed Ricci flow with surgery formula}) with initial value $u(0)=u^*$
exists for all time and converges to a discrete conformal factor with combinatorial $\alpha$-curvature $\overline{\mathbf{F}}$, which implies
the existence of a PH metric in the discrete conformal class $\mathcal{D}(d_0)$ with combinatorial $\alpha$-curvature $\overline{\mathbf{F}}$.

Set
\begin{equation*}
\mathbf{F}_{\alpha,\max}(u(0))=\max_{i\in V}\mathbf{F}_{\alpha,i}(u^\ast),\ M_{\alpha,\max}=\max_{i\in V}M_{\alpha,i}.
\end{equation*}
\begin{equation*}
\mathbf{F}_{\alpha,\min}(u(0))=\min_{i\in V}\mathbf{F}_{\alpha,i}(u^\ast),\  M_{\alpha,\min}=\min_{i\in V}M_{\alpha,i}.
\end{equation*}
Due to $\overline{\mathbf{F}}<\mathbf{F}_{\alpha,\min}(u(0))$, we have $M_{\alpha,i}(0)>0$ for all $v_i\in V$, which implies $M_{\alpha,i}(t)\geq 0$ for all $v_i\in V$ by Corollary \ref{application of maximum principle}.
Combining $A_i>0$, $M_{\alpha,i}(t)\geq 0$ and (\ref{the evolution of M alpha i}), we have
\begin{equation*}
\frac{dM_{\alpha,i}}{dt}\leq \sum_{j\sim i}\frac{B_{ij}}{w^\alpha_i}(M_{\alpha,j}-M_{\alpha,i})+\alpha M_{\alpha,i}(M_{\alpha,i}+\overline{\mathbf{F}}),
\end{equation*}
which implies
\begin{equation}\label{estimate M alpha 1}
0
\leq M_{\alpha,i}(t)
\leq \frac{\mathbf{\overline{F}}}{-1+(\frac{\mathbf{\overline{F}}}{M_{\alpha,\max}(0)}+1)e^{-\alpha \mathbf{\overline{F}}t}}
\leq \frac{\overline{\mathbf{F}}M_{\alpha,max}(u(0))}{\mathbf{F}_{\alpha,max}(u(0))}e^{\alpha \mathbf{\overline{F}}t}
\end{equation}
by Theorem \ref{maximum principle}.
Combining $\frac{du_i}{dt}=-M_{\alpha,i}$, $\alpha>0$ and $\overline{\mathbf{F}}<0$,
the estimation (\ref{estimate M alpha 1}) further implies $\{u(t)\}\subset\subset \mathbb{R}^N$.
Therefore, the solution of modified combinatorial $\alpha$-Yamabe flow with surgery (\ref{prescribed Ricci flow with surgery formula}) exists for all time and converges, which implies the existence of a PH metric in the conformal class of $\mathcal{D}(d_0)$ with combinatorial $\alpha$-curvature $\overline{\mathbf{F}}$.
\qed

\begin{remark}\label{remark on existence of constant alpha curvature metric}
Lemma \ref{Lemma 1} still holds under the condition that $\alpha<0$, $\overline{\mathbf{F}}$ is a positive constant and $0<\mathbf{F}_{\alpha,i}(u^\ast)<\overline{\mathbf{F}}$ for all $v_i\in V$,
the proof of which is the same as that for Lemma \ref{Lemma 1}.
\end{remark}
To handle more general cases, we need the following technique lemma.

\begin{lemma}\label{boundedness of u_i}
Suppose $\triangle v_iv_jv_k\in F$ is a hyperbolic triangle with discrete hyperbolic metrics $d_{ij},d_{jk},d_{ki}$.
If the discrete conformal factors $u_i,u_j,u_k$ for the hyperbolic triangle are uniformly bounded from below, i.e. $\exists C$ such that
$u_i, u_j, u_k\in [C, +\infty)$.
Then for any $\epsilon>0$, there exists a number $\widetilde{C}=\widetilde{C}(d_{ij},d_{jk},d_{ki}, C)$ such that if $u_i>\widetilde{C}$,
the inner angle $\alpha^{jk}_i$ is smaller than $\epsilon$.
\end{lemma}
\proof
We use the trick from \cite{GX1} to prove the lemma.
Due to
$$\sinh\frac{l_{ij}}{2}=\sinh \frac{d_{ij}}{2}e^{u_i+u_j}\geq\sinh \frac{d_{ij}}{2}e^{u_i+C},$$
we have $\sinh\frac{l_{ij}}{2}\rightarrow +\infty,\ \sinh\frac{l_{ik}}{2}\rightarrow +\infty$ uniformly as $u_i\rightarrow +\infty$, where $l_{rs}=(u*d)_{rs}$, $\{r,s,t\}=\{i,j,k\}$.
By the hyperbolic cosine law,
\begin{equation*}
\begin{aligned}
\cos \alpha^{jk}_i
=&\frac{\cosh l_{ij}\cosh l_{ik}-\cosh l_{jk}}{\sinh l_{ij}\sinh l_{ik}}\\
=&\frac{\cosh(l_{ij}+l_{ik})+\cosh(l_{ij}-l_{ik})-2\cosh l_{jk}}{\cosh(l_{ij}+l_{ik})-\cosh(l_{ij}-l_{ik})}\\
=&\frac{1+\lambda-2\mu}{1-\lambda},
\end{aligned}
\end{equation*}
where $\lambda=\frac{\cosh(l_{ij}-l_{ik})}{\cosh(l_{ij}+l_{ik})}$ and $\mu=\frac{\cosh l_{jk}}{\cosh(l_{ij}+l_{ik})}$.
To prove $\alpha^{jk}_i\rightarrow 0$ uniformly,
we just need to prove $\lambda, \mu\rightarrow 0$ uniformly as $u_i\rightarrow +\infty$.
Note that
\begin{equation*}
\begin{aligned}
0<&\lambda<\frac{\cosh l_{ij}}{\cosh(l_{ij}+l_{ik})}<\frac{1}{\cosh l_{ik}}\leq\frac{1}{\sinh \frac{d_{ik}}{2}e^{u_i+C}},\\
0<&\mu
=\frac{1}{(1+2\sinh^2\frac{l_{ij}}{2})(1+2\sinh^2\frac{l_{ik}}{2})}+
\frac{\sinh^2\frac{d_{jk}}{2}}{\sinh^2\frac{d_{ij}}{2}\sinh^2\frac{d_{ik}}{2}e^{4u_i}},
\end{aligned}
\end{equation*}
we have $\lambda, \mu\rightarrow 0$ uniformly as $u_i\rightarrow +\infty$.
Therefore, $\alpha^{jk}_i\rightarrow 0$ uniformly as $u_i\rightarrow +\infty$.
\qed

As an application of Lemma \ref{boundedness of u_i}, we have the following result.
\begin{lemma}\label{Lemma 2}
Suppose $(S,V)$ is a connected closed marked surface with a PH metric $d_0$, $\alpha>0$ is a constant and $\overline{\mathbf{F}}$ is a non-positive constant.
If there exists $u^*\in \mathcal{U}(d_0)$ such that $\mathbf{F}_{\alpha}(u^\ast)$ is a negative constant,
then there exists a PH metric in the discrete conformal class $\mathcal{D}(d_0)$
with combinatorial $\alpha$-curvature $\overline{\mathbf{F}}$.
\end{lemma}
\proof
Because $\mathbf{F}_{\alpha}(u^*)$ and $\overline{\mathbf{F}}$ are both constants, then
$$\overline{\mathbf{F}}=\mathbf{F}_{\alpha}(u^*)<0\quad \text{or} \quad
\overline{\mathbf{F}}<\mathbf{F}_{\alpha}(u^*)<0 \quad \text{or} \quad
\mathbf{F}_{\alpha}(u^*)<\overline{\mathbf{F}}\leq 0.$$

In the case that $\overline{\mathbf{F}}=\mathbf{F}_{\alpha}(u^*)<0$, $\mathbf{A}^{-1}(p,w(-2u^\ast))\in \mathcal{D}(d_0)$ is a PH metric
with combinatorial $\alpha$-curvature $\overline{\mathbf{F}}$.

In the case that $\overline{\mathbf{F}}<\mathbf{F}_{\alpha}(u^*)<0$, the conclusion follows from Lemma \ref{Lemma 1}.

In the case that $\mathbf{F}_{\alpha}(u^*)<\overline{\mathbf{F}}\leq0$, we
take $u^*$ as the initial value of the modified combinatorial $\alpha$-Yamabe flow with surgery (\ref{prescribed Ricci flow with surgery formula}).
Then $M_{\alpha}(0)<0$ is a constant.
By Corollary \ref{application of maximum principle}, $M_{\alpha,i}(t)\leq 0$ for all $v_i\in V$, which implies $\mathbf{F}_{\alpha,i}(u(t))\leq\overline{\mathbf{F}}\leq0$ for all $v_i\in V$.
Combining $A_i>0$, $M_{\alpha,i}(t)\leq 0$ and (\ref{the evolution of M alpha i}), we have
\begin{equation*}
\frac{dM_{\alpha,i}}{dt}\geq \sum_{j\sim i}\frac{B_{ij}}{w^\alpha_i}(M_{\alpha,j}-M_{\alpha,i})+\alpha M_{\alpha,i}(M_{\alpha,i}+\overline{\mathbf{F}}).
\end{equation*}
Therefore, by Theorem \ref{maximum principle},
\begin{equation}\label{estimate M alpha 2}
\frac{\mathbf{\overline{F}}}{-1+(\frac{\mathbf{\overline{F}}}{M_{\alpha}(0)}+1)e^{-\alpha \mathbf{\overline{F}}t}} \leq M_{\alpha,i}(t) \leq 0,\ \text{if}\ \mathbf{\overline{F}}<0,
\end{equation}
$$\frac{M_{\alpha}(0)}{-\alpha M_{\alpha}(0)t+1}\leq M_{\alpha,i}(t)\leq 0,\ \text{if}\ \mathbf{\overline{F}}=0.$$
In the case that $\mathbf{\overline{F}}<0$, by the the estimation (\ref{estimate M alpha 2}) and $\frac{du_i}{dt}=-M_{\alpha,i}(t)$,
we have $\{u(t)\}\subset\subset \mathbb{R}^N$.
In the case that $\mathbf{\overline{F}}=0$, we have $\frac{du_i}{dt}=-M_{\alpha,i}(t)\geq0$, which implies
$u_i(t)$ is uniformly bounded from below for all $v_i\in V$.
We claim that $u(t)$ is bounded from above. If not, then there exists at least one vertex $v_i\in V$ such that $\lim_{t\rightarrow +\infty}u_i(t)=+\infty$.
For the vertex $v_i$, by Lemma \ref{boundedness of u_i}, there exists $C_i>0$ such that whenever $u_i(t)>C_i$, the inner angle $\alpha^{jk}_i$ is smaller than $\frac{\pi}{N_i}$, where $N_i$ is the degree at the vertex $i$. Thus  $K_i=2\pi-\sum_{\triangle v_iv_jv_k\in F}\alpha^{jk}_i>\pi>0$ and then $\mathbf{F}_{\alpha,i}(u(t))>0$, which contradicts $\mathbf{F}_{\alpha,i}(u(t))\leq0$. Therefore, $u_i(t)$ is bounded from above.
 This completes the proof of the claim.
As a consequence, the solution of modified combinatorial $\alpha$-Yamabe flow with surgery (\ref{prescribed Ricci flow with surgery formula}) exists for all time and converges, which implies the existence of a PH metric in the discrete conformal class $\mathcal{D}(d_0)$ with combinatorial $\alpha$-curvature $\overline{\mathbf{F}}$.
\qed

\begin{remark}\label{Lemma 3,4}
Lemma \ref{Lemma 2} still holds under the condition that $\alpha<0$, $\overline{\mathbf{F}}$ is a positive constant and $\mathbf{F}_{\alpha}(u^\ast)$ is a positive constant. As we do not have a result paralleling to Lemma \ref{boundedness of u_i} under the condition that $u_i$ is bounded from above for all $v_i\in V$, the result in Lemma \ref{Lemma 2} can not be generalized to the case that $\alpha<0$, $\overline{\mathbf{F}}$ is a nonnegative constant and $\mathbf{F}_{\alpha}(u^\ast)$ is a positive constant.
\end{remark}

As a direct corollary of Lemma \ref{Lemma 1}, Lemma \ref{Lemma 2}, Remark \ref{remark on existence of constant alpha curvature metric}
and Remark \ref{Lemma 3,4}, we have the following result.
\begin{corollary}\label{corollary}
Suppose $(S,V)$ is a connected closed marked surface with a PH metric $d_0$ and $\alpha$ is a constant.
\begin{description}
  \item[(1)] If $\alpha>0$ and there exists $u^*\in \mathcal{U}(d_0)$ such that $\mathbf{F}_{\alpha,i}(u^\ast)<0$ for all $v_i\in V$,
  then for any nonpositive constant $\overline{\mathbf{F}}$,
  there exists a PH metric in $\mathcal{D}(d_0)$ with combinatorial $\alpha$-curvature $\overline{\mathbf{F}}$;
  \item[(2)] If $\alpha<0$ and there exists $u^*\in \mathcal{U}(d_0)$ such that $\mathbf{F}_{\alpha,i}(u^\ast)>0$ for all $v_i\in V$,
  then for any positive constant $\overline{\mathbf{F}}$,
  there exists a PH metric in $\mathcal{D}(d_0)$ with combinatorial $\alpha$-curvature $\overline{\mathbf{F}}$.
\end{description}
\end{corollary}

The following result shows that the constant condition on $\overline{\mathbf{F}}$ in Corollary \ref{corollary} could be removed.
\begin{theorem}\label{main theorem}
Suppose $(S,V)$ is a connected closed marked surface with a PH metric $d_0$, $\alpha$ is a constant and $\overline{\mathbf{F}}: V\rightarrow \mathbb{R}$ is a function defined on $V$. If there exists $u^\ast\in \mathcal{U}(d_0)$ such that one of the following conditions is satisfied :
\begin{description}
  \item[(1)] $\alpha>0,\ \overline{\mathbf{F}}_i\leq0,\ \mathbf{F}_{\alpha,i}(u^\ast)<0$ for all $v_i\in V$,
  \item[(2)] $\alpha<0,\ \overline{\mathbf{F}}_i>0,\ \mathbf{F}_{\alpha,i}(u^\ast)>0$ for all $v_i\in V$,
\end{description}
then there exists a PH metric in the discrete conformal class $\mathcal{D}(d_0)$ with combinatorial $\alpha$-curvature $\overline{\mathbf{F}}$.
\end{theorem}
\proof
In the case that $\alpha>0,\ \overline{\mathbf{F}}\leq0,\ \mathbf{F}_{\alpha,i}(u^\ast)<0$, set
$$\overline{\mathbf{F}}_{\max}=\max_{v_i\in V}\overline{\mathbf{F}}_i,\ \overline{\mathbf{F}}_{\min}=\min_{v_i\in V}\overline{\mathbf{F}}_i.$$
By Corollary \ref{corollary}, there exists $\overline{u}\in \mathcal{U}(d_0)$
with negative constant combinatorial $\alpha$-curvature $\mathbf{F}_{\alpha}(\overline{u})<\overline{\mathbf{F}}_{\min}<0$.
Taking $\overline{u}$ as the initial value of the modified combinatorial $\alpha$-Yamabe flow with surgery (\ref{prescribed Ricci flow with surgery formula}),
then $M_{\alpha,i}(0)<0$ for all $v_i\in V$, which implies $M_{\alpha,i}(t)\leq 0$ for all $v_i\in V$ by Corollary \ref{application of maximum principle}.
Combining $\alpha>0$, $\overline{\mathbf{F}}\leq0$, $A_i>0$, $M_{\alpha,i}(t)\leq 0$ and (\ref{the evolution of M alpha i}), we have
\begin{equation*}
\begin{aligned}
\frac{dM_{\alpha,i}}{dt}
&\geq \sum_{j\sim i}\frac{B_{ij}}{w^\alpha_i}(M_{\alpha,j}-M_{\alpha,i})+\alpha M_{\alpha,i}(M_{\alpha,i}+\overline{\mathbf{F}}_i) \\
&\geq \sum_{j\sim i}\frac{B_{ij}}{w^\alpha_i}(M_{\alpha,j}-M_{\alpha,i})+\alpha M_{\alpha,i}(M_{\alpha,i}+\overline{\mathbf{F}}_{\max}).
\end{aligned}
\end{equation*}
Therefore, by Theorem \ref{maximum principle},
$$\frac{\mathbf{\overline{F}}_{\max}}{-1+(\frac{\mathbf{\overline{F}}_{\max}}{M_{\alpha,\min}(0)}+1)e^{-\alpha \mathbf{\overline{F}}_{\max}t}} \leq M_{\alpha,i}(t) \leq 0,\ \text{if}\ \overline{\mathbf{F}}_{\max}<0,$$
$$\frac{M_{\alpha,\min}(0)}{-\alpha M_{\alpha,\min}(0)t+1}\leq M_{\alpha,i}(t)\leq 0,\ \text{if}\  \overline{\mathbf{F}}_{\max}= 0.$$
Following the arguments in the proof of Lemma \ref{Lemma 2},
the solution of the modified combinatorial $\alpha$-Yamabe flow with surgery (\ref{prescribed Ricci flow with surgery formula}) stays in a compact subset of $\mathbb{R}^N$.
Therefore, the solution of (\ref{prescribed Ricci flow with surgery formula}) exists for all time and converges,
which implies the existence of a PH metric in the discrete conformal class of $\mathcal{D}(d_0)$
with combinatorial $\alpha$-curvature $\overline{\mathbf{F}}$.

The same arguments apply to the case that $\alpha<0,\ \overline{\mathbf{F}}>0,\ \mathbf{F}_{\alpha,i}(u(0))>0$. We omit the details here.
\qed

\textbf{Proof of Theorem \ref{discrete uniformization theorem}:}
In the case of $\alpha>0,\ \chi(S)<0,\ \overline{\mathbf{F}}\leq0$,
  by Theorem \ref{Gu. discrete uniformization theorem}, there exists $\overline{u}\in \mathcal{U}(d_0)$ such that $\sum^N_i K_i^\ast(\overline{u})>2\pi \chi(S)$ and $K_i^\ast(\overline{u})<0$ for all $v_i\in V$. Then $\mathbf{F}_{\alpha,i}(u^\ast)<0$ for all $v_i\in V$.
By Theorem \ref{main theorem}, there exists a PH metric in the discrete conformal class $\mathcal{D}(d_0)$ with
non-positive combinatorial $\alpha$-curvature $\overline{\mathbf{F}}$.

In the case of $\alpha<0, \overline{\mathbf{F}}>0$,
by Theorem \ref{Gu. discrete uniformization theorem},
there exists $u^*\in \mathcal{U}(d_0)$ such that $\mathbf{F}_{\alpha,i}(u^\ast)>0$ for all $v_i\in V$.
By Theorem \ref{main theorem}, there exists a PH metric in the discrete conformal class $\mathcal{D}(d_0)$
with positive combinatorial $\alpha$-curvature $\overline{\mathbf{F}}$.

The case of $\alpha=0$ follows from Theorem \ref{Gu. discrete uniformization theorem}.
\qed

As a corollary of Theorem \ref{equivalent theorem 1} and Theorem \ref{discrete uniformization theorem},
we have the following longtime existence and convergence for the solution of modified
combinatorial $\alpha$-Yamabe flow (\ref{prescribed Ricci flow with surgery formula}),
which is the first part of Theorem \ref{last theorem}.
\begin{theorem}
Suppose $(S,V)$ is a connected closed marked surface with a PH metric $d_0$, $\alpha$ is a constant and $\overline{\mathbf{F}}: V\rightarrow \mathbb{R}$ is a function defined on $V$.
If one of the conditions $\mathbf{(1)(2)(3)}$ in Theorem \ref{discrete uniformization theorem} is satisfied,
then the solution of modified combinatorial $\alpha$-Yamabe flow with surgery (\ref{prescribed Ricci flow with surgery formula}) exists for all time and converges exponentially fast to a PH metric with combinatorial $\alpha$-curvature $\overline{\mathbf{F}}$.
\end{theorem}

\subsection{Combinatorial $\alpha$-Calabi flow with surgery}
Analogy to Definition \ref{Ricci flow with surgery definition} for combinatorial $\alpha$-Yamabe flow with surgery,
we define the following combinatorial $\alpha$-Calabi flow with surgery.

\begin{definition}\label{Calabi flow with surgery definition}
Suppose $(S,V)$ is a connected closed marked surface with a PH metric $d_0$ and $\alpha\in \mathbb{R}$ is a constant. The combinatorial $\alpha$-Calabi flow with surgery is defined to be
\begin{eqnarray}\label{Calabi flow with surgery formula}
\begin{cases}
\frac{du_i}{dt}=\Delta_\alpha \mathbf{F}_{\alpha,i},\\
u_i(0)=0,
\end{cases}
\end{eqnarray}
where $\Delta_\alpha$ is the $\alpha$-Laplace operator defined by (\ref{alpha laplacian definition}).
\end{definition}
One can use the following modified combinatorial $\alpha$-Calabi flow with surgery
to study prescribed combinatorial $\alpha$-curvature problem on polyhedral surfaces
\begin{equation}\label{prescribed Calabi flow with surgery formula}
\frac{du_i}{dt}=\Delta_\alpha (\mathbf{F}_{\alpha}-\overline{\mathbf{F}})_i,
\end{equation}
where $\overline{\mathbf{F}}$ is a function defined on $V$.

\begin{theorem}\label{equivalent theorem 2}
Suppose $(S,V)$ is a connected closed marked surface with a PH metric $d_0$, $\alpha\in \mathbb{R}$ is a constant
and $\overline{\mathbf{F}}: V\rightarrow \mathbb{R}$ is a function defined on $V$ with $\alpha \overline{\mathbf{F}}\leq 0$.
Then there exists a PH metric in the discrete conformal class $\mathcal{D}(d_0)$ with combinatorial $\alpha$-curvature $\overline{\mathbf{F}}$ if and only if the solution of modified combinatorial $\alpha$-Calabi flow with surgery (\ref{prescribed Calabi flow with surgery formula})
exists for all time and converges exponentially fast to a PH metric with combinatorial $\alpha$-curvature $\overline{\mathbf{F}}$.
\end{theorem}
\proof
Suppose the solution $u(t)$ of the modified combinatorial $\alpha$-Calabi flow with surgery (\ref{prescribed Calabi flow with surgery formula})
converges to $u^{\ast}\in \mathcal{U}(d_0)$ as $t\rightarrow +\infty$, then
$\mathbf{F}_\alpha(u^*)=\lim_{t\rightarrow +\infty}\mathbf{F}_\alpha(u(t))$ by the $C^1$-smoothness of $\mathbf{F}$.
Furthermore, there exists a sequence $\xi_n\in(n,n+1)$ such that
\begin{equation*}
u_i(n+1)-u_i(n)=u'_i(\xi_n)=\Delta_\alpha(\mathbf{F}_{\alpha}(u(\xi_n))-\mathbf{\overline{F}})_i\rightarrow 0,\ \text{as}\ n\rightarrow +\infty,
\end{equation*}
which implies $\mathbf{F}_{\alpha,i}(u^*)=\lim_{n\rightarrow +\infty}\mathbf{F}_{\alpha,i}(u(\xi_n))=\overline{\mathbf{F}}_i$ for all $v_i\in V$ and $u^*\in \mathcal{U}(d_0)$ corresponds to a PH metric with combinatorial $\alpha$-curvature $\overline{\mathbf{F}}$.

Conversely, suppose there exists a PH metric in the conformal class $\mathcal{D}(d_0)$ with combinatorial  $\alpha$-curvature $\overline{\mathbf{F}}$
and $u^*$ is the corresponding discrete conformal factor in $\mathcal{U}(d_0)$. Set
\begin{equation*}
W_\alpha(u)=W(u)-\int_{u^*}^{u}\sum_{i=1}^{N}\overline{\mathbf{F}}_iw^{\alpha}_idu_i
=\int^u_{u^*}\sum_{i=1}^N\mathbf{F}_idu_i-\int_{u^*}^{u}\sum_{i=1}^{N}\overline{\mathbf{F}}_iw^{\alpha}_idu_i.
\end{equation*}
Then $W_\alpha(u)$ is a $C^2$-smooth strictly convex function defined on $\mathbb{R}^n$ by Proposition \ref{W} and $\alpha \overline{\mathbf{F}}\leq 0$.
Furthermore, $W_\alpha(u^*)=0,\ \nabla W_\alpha(u^*)=0$ and $\mathrm{Hess} W_\alpha>0$, which implies $W_\alpha(u)\geq W_\alpha(u^\ast)=0$ and $\lim_{u\rightarrow \infty}W_\alpha(u)=+\infty$. By direct calculations, we have
\begin{equation*}
\frac{dW_\alpha(u(t))}{dt}
=\sum^N_{i=1}\frac{\partial W_\alpha}{\partial u_i}\frac{du_i}{dt}
=\sum^N_{i=1}(\mathbf{F}_{\alpha}-\overline{\mathbf{F}})_iw^\alpha_i
\Delta_\alpha(\mathbf{F}_{\alpha}-\overline{\mathbf{F}})_i
=-(\mathbf{F}_\alpha-\overline{\mathbf{F}})^T\cdot \widetilde{L}\cdot (\mathbf{F}_\alpha-\overline{\mathbf{F}})\leq 0
\end{equation*}
by Proposition \ref{W},
which implies $0\leq W_\alpha(u(t))\leq W_\alpha(u(0))$. Thus $\{u(t)\}\subset\subset \mathbb{R}^N$  by $\lim_{u\rightarrow \infty}W_\alpha(u)=+\infty$,
which implies the solution of the modified combinatorial $\alpha$-Calabi flow with surgery (\ref{prescribed Calabi flow with surgery formula}) exists for all time and $W_\alpha(u(t))$ converges.
Moreover, there exists a sequence $\xi_n\in(n,n+1)$ such that as $n\rightarrow +\infty$,
\begin{equation*}
\begin{aligned}
&W_\alpha(u(n+1))-W_\alpha(u(n))=(W_\alpha(u(t))'|_{\xi_n}=\nabla W_\alpha\cdot\frac{du_i}{dt}|_{\xi_n}\\
=&\sum^N_{i=1}(\mathbf{F}_{\alpha}-\overline{\mathbf{F}})_iw^\alpha_i
\Delta_\alpha(\mathbf{F}_{\alpha}-\overline{\mathbf{F}})_i|_{\xi_n}
=-(\mathbf{F}_\alpha-\overline{\mathbf{F}})^T\cdot \widetilde{L}\cdot (\mathbf{F}_\alpha-\overline{\mathbf{F}})|_{\xi_n}\rightarrow 0.
\end{aligned}
\end{equation*}
Then $\lim_{n\rightarrow +\infty}\mathbf{F}_{\alpha,i}(u(\xi_n))=\overline{\mathbf{F}}_i=\mathbf{F}_{\alpha,i}(u^*)$ for all $v_i\in V$, which implies $\lim_{n\rightarrow +\infty}u(\xi_n)=u^*$ by Theorem \ref{F global rigidity theorem}.

Similar to the proof of Theorem \ref{alpha flows with small initial energy}, set $\Gamma(u)=\Delta_\alpha (\mathbf{F}_{\alpha}-\overline{\mathbf{F}})$,
then $D\Gamma|_{u=u^\ast}$ has $N$ negative eigenvalues, which implies that $u^*$ is a local attractor of (\ref{prescribed Calabi flow with surgery formula}). Then the conclusion follows by Lyapunov Stability Theorem (\cite{Pontryagin}, Chapter 5).
\qed

As a corollary of Theorem \ref{equivalent theorem 2} and Theorem \ref{discrete uniformization theorem}, we have the following
longtime existence and convergence for the solution of modified combinatorial $\alpha$-Calabi flow with surgery (\ref{prescribed Calabi flow with surgery formula}),
which is the second part of Theorem \ref{last theorem}.
\begin{theorem}
Suppose $(S,V)$ is a connected closed marked surface with a PH metric $d_0$,
$\alpha$ is a constant and $\overline{\mathbf{F}}: V\rightarrow \mathbb{R}$ is a function defined on $V$.
If one of the conditions $\mathbf{(1)(2)(3)}$ in Theorem \ref{discrete uniformization theorem} is satisfied, then
the solution of modified combinatorial $\alpha$-Calabi flow with surgery (\ref{prescribed Calabi flow with surgery formula})
exists for all time and converges exponentially fast to a PH metric in $\mathcal{D}(d_0)$
with  combinatorial $\alpha$-curvature $\overline{\mathbf{F}}$.
\end{theorem}

\appendix
\section{Proof of Lemma \ref{L=A+L_B}}\label{appendix a}
We use $j\sim i$ to denote that the vertices $v_i$ and $v_j$ are adjacent.
By Lemma \ref{constants to be continuous}, $\frac{\partial \alpha^{jk}_i}{\partial u_j}=\frac{\partial \alpha^{ik}_j}{\partial u_i}$ and $\frac{\partial K_i}{\partial u_j}=\frac{\partial K_j}{\partial u_i}$.
\begin{description}
  \item[(1)] If $j\nsim i$ and $j\neq i$, then $\frac{\partial K_i}{\partial u_j}=0.$
  \item[(2)] If $j\sim i$, then
\begin{equation*}
\frac{\partial K_i}{\partial u_j}
=\frac{\partial (2\pi-\sum_{\triangle v_iv_jv_k\in F}\alpha^{jk}_i)}{\partial u_j}
=-\sum_{\triangle v_iv_jv_k\in F}\frac{\partial \alpha^{jk}_i}{\partial u_j}
=-(\frac{\partial \alpha^{jk}_i}{\partial u_j}+\frac{\partial \alpha^{jl}_i}{\partial u_j}),
\end{equation*}
where $\triangle v_iv_jv_k$ and $\triangle v_iv_jv_l$ are adjacent triangles sharing a common edge $\{v_iv_j\}$.
Set $B_{ij}=\frac{\partial \alpha^{jk}_i}{\partial u_j}+\frac{\partial \alpha^{jl}_i}{\partial u_j}$, then
$\frac{\partial K_i}{\partial u_j}=-B_{ij}.$
  \item[(3)] If $j=i$, then
\begin{equation*}
\begin{aligned}
\frac{\partial K_i}{\partial u_i}
=&\sum_{\triangle v_iv_jv_k\in F}\frac{\partial \left(\alpha^{ik}_j+\alpha^{ij}_k+\text{Area}(\triangle v_iv_jv_k)\right)}{\partial u_i}\\
=&\sum_{\triangle v_iv_jv_k\in F}\left(\frac{\partial \alpha^{jk}_i}{\partial u_j}+\frac{\partial \alpha^{jk}_i}{\partial u_k}\right)
+\sum_{\triangle v_iv_jv_k\in F}\frac{\partial \text{Area}(\triangle v_iv_jv_k)}{\partial u_i}\\
=&\sum_{j\sim i}\left(\frac{\partial \alpha^{jk}_i}{\partial u_j}+\frac{\partial \alpha^{jl}_i}{\partial u_j}\right)+\frac{\partial}{\partial u_i}\left(\sum_{\triangle v_iv_jv_k\in F}\text{Area}(\triangle v_iv_jv_k)\right).
\end{aligned}
\end{equation*}
\end{description}
Set $A_i=\frac{\partial}{\partial u_i}\left(\sum_{\triangle v_iv_jv_k\in F}\text{Area}(\triangle v_iv_jv_k)\right)$, then
\begin{equation}\label{(3)}
\frac{\partial K_i}{\partial u_i}=A_i+\sum_{j\sim i}B_{ij}.
\end{equation}
Set $L_A=\text{diag}\{A_1,...,A_N\}$ and $L_B=((L_B)_{ij})_{N\times N}$, where
\begin{eqnarray*}
(L_B)_{ij}=
\begin{cases}
\sum_{k\sim i}B_{ik},  &{j=i},\\
-B_{ij}, &{j\sim i},\\
0, &{j\nsim i, j\neq i}.
\end{cases}
\end{eqnarray*}
Then $L=L_A+L_B$.

To get the expression of $B_{ij}$ in (\ref{B}), we need the following lemma.
\begin{lemma}\label{2}
Suppose $\triangle v_iv_jv_k$ is a hyperbolic triangle with a PH metric $d_0$ and $u$ is a discrete conformal factor for the
triangle. Denote $d_{st}:=(u*d_0)_{st}$ as the edge length of $\{st\}$ in the triangle  $\triangle v_iv_jv_k$
and $\alpha^{jk}_i, \alpha^{ik}_j, \alpha^{ij}_k$ as the inner angles facing $\{jk\},\{ik\},\{ij\}$ respectively. Then
\begin{equation*}
\frac{\partial \alpha^{jk}_i}{\partial u_j}
=\frac{1}{\cosh^2 \frac{d_{ij}}{2}} \tan \frac{\alpha^{jk}_i+\alpha^{ik}_j-\alpha^{ij}_k}{2}.
\end{equation*}
\end{lemma}
\proof
By the derivative cosine law (see Lemma A1 in \cite{Chow-Luo} for example),
\begin{equation*}
\begin{aligned}
\frac{\partial \alpha^{jk}_i}{\partial d_{jk}}=\frac{\sinh d_{jk}}{A},\
\frac{\partial \alpha^{jk}_i}{\partial d_{ik}}=\frac{-\sinh d_{jk}\cos\alpha^{ij}_k}{A},\
\frac{\partial \alpha^{jk}_i}{\partial d_{ij}}=\frac{-\sinh d_{jk}\cos\alpha^{ik}_j}{A},
\end{aligned}
\end{equation*}
where $A=\sinh d_{ik}\sinh d_{ij}\sin\alpha^{jk}_i$.
By (\ref{discrete conformal equivalent formula}), we have
$$\frac{\partial d_{jk}}{\partial u_i}=0,\ \frac{\partial d_{jk}}{\partial u_j}=\frac{\partial d_{jk}}{\partial u_k}=2\tanh\frac{d_{jk}}{2}.$$
According to the chain rules,
\begin{equation}\label{equation in proof apendix 2}
\begin{aligned}
\frac{\partial \alpha^{jk}_i}{\partial u_j}
=&\frac{\partial \alpha^{jk}_i}{\partial d_{jk}}\frac{\partial d_{jk}}{\partial u_j}+\frac{\partial \alpha^{jk}_i}{\partial d_{ik}}\frac{\partial d_{ik}}{\partial u_j}+\frac{\partial \alpha^{jk}_i}{\partial d_{ij}}\frac{\partial d_{ij}}{\partial u_j}\\
=&\frac{2\sinh d_{jk}}{A}\tanh\frac{d_{jk}}{2}-\frac{2\sinh d_{jk}\cos\alpha^{ik}_j}{A}\tanh\frac{d_{ij}}{2}\\
=&\frac{2(\cosh d_{jk}+\cosh d_{ik}-\cosh d_{ij}-1)}{A(\cosh d_{ij}+1)}.
\end{aligned}
\end{equation}
By the hyperbolic cosine law expressing $d_{ik}$ in terms of $\alpha^{jk}_i, \alpha^{ik}_j, \alpha^{ij}_k$, we have
\begin{equation}\label{equation in proof apendix 3}
\begin{aligned}
\cosh^2 \frac{d_{ik}}{2}
=&\frac{1}{2}(\cosh d_{ik}+1)\\
=&\frac{1}{2}\left(\frac{\cos \alpha^{ik}_j+\cos \alpha^{ij}_k\cos \alpha^{jk}_i}{\sin \alpha^{ij}_k \sin \alpha^{jk}_i}+1\right)\\
=&\frac{\cos\frac{\alpha^{ik}_j+\alpha^{ij}_k-\alpha^{jk}_i}{2}\cos\frac{\alpha^{ik}_j-\alpha^{ij}_k+
\alpha^{jk}_i}{2}}{\sin \alpha^{ij}_k \sin \alpha^{jk}_i},
\end{aligned}
\end{equation}
which implies
\begin{equation}\label{equation in proof apendix}
\begin{aligned}
&\cosh\frac{d_{jk}}{2}\cosh\frac{d_{ik}}{2}\sin \alpha^{ij}_k\\
=&\sqrt{\frac{\cos\frac{\alpha^{jk}_i+\alpha^{ik}_j-\alpha^{ij}_k}{2}\cos\frac{\alpha^{jk}_i-\alpha^{ik}_j+
\alpha^{ij}_k}{2}}{\sin \alpha^{ik}_j \sin \alpha^{ij}_k}}\cdot \sqrt{\frac{\cos\frac{\alpha^{ik}_j+\alpha^{ij}_k-\alpha^{jk}_i}{2}\cos\frac{\alpha^{ik}_j-\alpha^{ij}_k+
\alpha^{jk}_i}{2}}{\sin \alpha^{ij}_k \sin \alpha^{jk}_i}}\sin \alpha^{ij}_k\\
=&\cosh\frac{d_{ij}}{2}\cos\frac{\alpha^{jk}_i+\alpha^{ik}_j-\alpha^{ij}_k}{2},
\end{aligned}
\end{equation}
where (\ref{equation in proof apendix 3}) is used in the last line.
Recall the following formula obtained by Gu-Guo-Luo-Sun-Wu (\cite{Gu2} Lemma 11)
\begin{equation}\label{GGLSW's formula}
2\sin\frac{\alpha^{jk}_i+\alpha^{ik}_j-\alpha^{ij}_k}{2}\cdot \cosh \frac{d_{ij}}{2}=
\frac{\sinh^2 \frac{d_{jk}}{2}+\sinh^2 \frac{d_{ik}}{2}-\sinh^2 \frac{d_{ij}}{2}}{\sinh \frac{d_{jk}}{2}\sinh \frac{d_{ik}}{2}}.
\end{equation}
Then (\ref{equation in proof apendix 2}) implies
\begin{equation*}
\begin{aligned}
\frac{\partial \alpha^{jk}_i}{\partial u_j}
=&\frac{2(\cosh d_{jk}+\cosh d_{ik}-\cosh d_{ij}-1)}{\sinh d_{ik}\sinh d_{jk}\sin\alpha^{ij}_k(\cosh d_{ij}+1)}\\
=&\frac{4(\sinh^2 \frac{d_{jk}}{2}+\sinh^2 \frac{d_{ik}}{2}-\sinh^2 \frac{d_{ij}}{2})}{\sinh d_{ik}\sinh d_{jk}\sin\alpha^{ij}_k(\cosh d_{ij}+1)}\\
=&\frac{8\sinh \frac{d_{jk}}{2}\sinh \frac{d_{ik}}{2}\cosh \frac{d_{ij}}{2}}{\sinh d_{ik}\sinh d_{jk}\sin\alpha^{ij}_k(\cosh d_{ij}+1)}\sin\frac{\alpha^{jk}_i+\alpha^{ik}_j-\alpha^{ij}_k}{2}\\
=&\frac{1}{\cosh\frac{d_{ij}}{2}\cosh\frac{d_{jk}}{2}\cosh\frac{d_{ik}}{2}\sin \alpha^{ij}_k}\sin\frac{\alpha^{jk}_i+\alpha^{ik}_j-\alpha^{ij}_k}{2}\\
=&\frac{1}{\cosh^2 \frac{d_{ij}}{2}} \tan \frac{\alpha^{jk}_i+\alpha^{ik}_j-\alpha^{ij}_k}{2},
\end{aligned}
\end{equation*}
where (\ref{GGLSW's formula}) is used in the third line and (\ref{equation in proof apendix}) is used in the last line.
\qed

\begin{lemma}
Suppose $(S,V,\mathcal{T})$ is a connected closed triangulated surface with a PH metric $d_0$.
$d:=u*d_0$ is a PH metric defined by a function $u:V\rightarrow \mathbb{R}$.
Then the triangulation $\mathcal{T}$ is Delaunay in $d$ if and only if $B_{ij}\geq 0$ for any $\{v_iv_j\}\in E$,
which implies $L_B$ is positive semi-definite.
Furthermore,
\begin{equation}\label{relation of B_ij and A_i}
\begin{aligned}
A_i=\sum_{j\sim i}B_{ij}(\cosh d_{ij}-1),
\end{aligned}
\end{equation}
which is nonnegative if  the triangulation $\mathcal{T}$ is Delaunay in $d$.
\end{lemma}
\proof
By Lemma \ref{2}, for $j\sim i$,
\begin{equation}\label{B_ij}
\begin{aligned}
B_{ij}
=&\frac{\partial \alpha^{jk}_i}{\partial u_j}+\frac{\partial \alpha^{jl}_i}{\partial u_j}\\
=&\frac{1}{\cosh^2 \frac{d_{ij}}{2}}\left(\tan \frac{\alpha^{jk}_i+\alpha^{ik}_j-\alpha^{ij}_k}{2}+\tan \frac{\alpha^{jl}_i+\alpha^{il}_j-\alpha^{ij}_l}{2}\right)\\
=&\frac{1}{\cosh^2 \frac{d_{ij}}{2}}\cdot
\frac{\sin
\frac{\alpha^{jk}_i+\alpha^{ik}_j-\alpha^{ij}_k+\alpha^{jl}_i+\alpha^{il}_j-\alpha^{ij}_l}{2}}
{\sin\frac{\pi-\alpha^{jk}_i-\alpha^{ik}_j+\alpha^{ij}_k}{2}
\sin\frac{\pi-\alpha^{jl}_i-\alpha^{il}_j+\alpha^{ij}_l}{2}}.
\end{aligned}
\end{equation}
Note that the triangulation $\mathcal{T}$ is Delaunay if and only if $\alpha^{jk}_i+\alpha^{ik}_j+\alpha^{jl}_i+\alpha^{il}_j\geq \alpha^{ij}_l+\alpha^{ij}_k$ \cite{Gu2,GL}, which is equivalent to $B_{ij}\geq 0$.

For any $x\in \mathbb{R}^N$, we have
$$x^TL_Bx=\sum_{i,j=1}^N(L_B)_{ij}x_ix_j=\sum_{i\sim j}(-B_{ij}x_ix_j)+\sum_{i\sim j}x_i^2B_{ij}=\frac{1}{2}\sum_{i\sim j}B_{ij}(x_i-x_j)^2\geq 0,$$
which implies $L_B$ is positive semi-definite.
One can also refer to Lemma 3.10 in \cite{Chow-Luo} for the positive semi-definiteness of $L_B$.

Recall the following formula obtained by Glickenstein-Thomas (\cite{GT} Proposition 9)
\begin{equation*}
\frac{\partial \text{Area}(\triangle v_iv_jv_k)}{\partial u_k}=\frac{\partial \alpha^{jk}_i}{\partial u_k}(\cosh d_{ik}-1)+\frac{\partial \alpha^{ik}_j}{\partial u_k}(\cosh d_{jk}-1).
\end{equation*}
Then
\begin{equation*}
\begin{aligned}
A_i
=&\sum_{\triangle v_iv_jv_k\in F}\frac{\partial \text{Area}(\triangle v_iv_jv_k)}{\partial u_i}\\
=&\sum_{\triangle v_iv_jv_k\in F}\frac{\partial \alpha^{jk}_i}{\partial u_j}(\cosh d_{ij}-1)+
\sum_{\triangle v_iv_jv_k\in F}\frac{\partial \alpha^{jk}_i}{\partial u_k}(\cosh d_{ik}-1)\\
=&\sum_{j\sim i}\left(\frac{\partial \alpha^{jk}_i}{\partial u_j}+\frac{\partial \alpha^{jl}_i}{\partial u_j}\right)(\cosh d_{ij}-1)\\
=&\sum_{j\sim i}B_{ij}(\cosh d_{ij}-1).
\end{aligned}
\end{equation*}
Therefore, $B_{ij}\geq 0$ implies $A_i\geq 0$.
\qed

\begin{remark}
Submitting (\ref{B}) into (\ref{relation of B_ij and A_i}) gives
\begin{equation}\label{A_i in appendix}
\begin{aligned}
A_i
=2\sum_{j\sim i}\tanh^2\frac{d_{ij}}{2}\cdot\left(\tan \frac{\alpha^{jk}_i+\alpha^{ik}_j-\alpha^{ij}_k}{2}+\tan \frac{\alpha^{jl}_i+\alpha^{il}_j-\alpha^{ij}_l}{2}\right).
\end{aligned}
\end{equation}
Then the decomposition $L=L_A+L_B$ with $A_i$ given by (\ref{A_i in appendix}) and $B_{ij}$ given by (\ref{B}) is equivalent to the formula obtained by
Bobenko-Pinkall-Springborn (\cite{BPS}, Proposition 6.1.7).
\end{remark}

We further have the following result on the sign of $A_i$ under the Delaunay condition, which completes the proof of Lemma \ref{L=A+L_B}.

\begin{lemma}
Suppose $(S,V,\mathcal{T})$ is a connected closed triangulated surface with a PH metric $d_0$.
$d:=u*d_0$ is a PH metric defined by a function $u:V\rightarrow \mathbb{R}$.
If the triangulation $\mathcal{T}$ is Delaunay in $d$, then
$A_i>0$ for all $v_i\in V$, which implies $L_A$ is positive definite.
\end{lemma}
\proof
Combining (\ref{(3)}) and (\ref{relation of B_ij and A_i}) gives
\begin{equation}\label{3}
\frac{\partial K_i}{\partial u_i}=A_i+\sum_{j\sim i}B_{ij}=\sum_{j\sim i}B_{ij}\cosh d_{ij}.
\end{equation}
For a single triangle $\triangle v_iv_jv_k$, by direct calculations, we have
\begin{equation}\label{derivative of alpha i to ui}
\frac{\partial \alpha^{jk}_i}{\partial u_i}
=\frac{\cosh^2 d_{ik}+\cosh^2 d_{ij}-2\cosh d_{jk}\cosh d_{ik}\cosh d_{ij}+(1-\cosh d_{jk})(\cosh d_{ik}+\cosh d_{ij})}{A(1+\cosh d_{ik})(1+\cosh d_{ij})},
\end{equation}
where $A=\sinh d_{ik}\sinh d_{ij}\sin\alpha^{jk}_i$. Note that
\begin{equation*}
\begin{aligned}
\sin^2 \alpha^{jk}_i
=&1-\left(\frac{\cosh d_{ik}\cosh d_{ij}-\cosh d_{jk}}{\sinh d_{ik}\sinh d_{ij}}\right)^2\\
=&\frac{1}{\sinh^2 d_{ik}\sinh^2 d_{ij}}(-\cosh^2 d_{jk}-\cosh^2 d_{ik}-\cosh^2 d_{ij}+1+2\cosh d_{jk}\cosh d_{ik}\cosh d_{ij}),
\end{aligned}
\end{equation*}
which implies
\begin{equation*}
\begin{aligned}
&-\cosh^2 d_{ik}-\cosh^2 d_{ij}+2\cosh d_{jk}\cosh d_{ik}\cosh d_{ij}+(\cosh d_{jk}-1)(\cosh d_{ik}+\cosh d_{ij})\\
=&\sinh^2 d_{ik}\sinh^2 d_{ij}\sin^2 \alpha^{jk}_i+\sinh^2 d_{jk}+(\cosh d_{jk}-1)(\cosh d_{ik}+\cosh d_{ij})\\
>&0,
\end{aligned}
\end{equation*}
we have $\frac{\partial \alpha^{jk}_i}{\partial u_i}< 0$ by (\ref{derivative of alpha i to ui}).
Then
\begin{equation*}
\frac{\partial K_i}{\partial u_i}
=\frac{\partial (2\pi-\sum_{\triangle v_iv_jv_k\in F}\alpha^{jk}_i)}{\partial u_i}
=-\sum_{\triangle v_iv_jv_k\in F}\frac{\partial \alpha^{jk}_i}{\partial u_i}>0,
\end{equation*}
which implies $\sum_{j\sim i}B_{ij}>0$ by (\ref{3}) and $B_{ij}\geq 0$.
Therefore, $A_i=\sum_{j\sim i}B_{ij}(\cosh d_{ij}-1)>0$ by (\ref{relation of B_ij and A_i}).
\qed

After the paper was finished, Dr. Tianqi Wu and Dr. Xiaoping Zhu told us that there was a similar proof of Lemma \ref{L=A+L_B}
in  \cite{WZ}. We thank Dr. Tianqi Wu and Dr. Xiaoping Zhu for informing us this and communications on related topics.
However, our proof is different. So we present our proof here.

\bigskip

(Xu Xu) School of Mathematics and Statistics, Wuhan University, Wuhan 430072, P.R. China

E-mail: xuxu2@whu.edu.cn\\[2pt]

(Chao Zheng) School of Mathematics and Statistics, Wuhan University, Wuhan 430072, P.R. China

E-mail: 2019202010023@whu.edu.cn\\[2pt]

\end{document}